\documentclass[11pt]{article}
\usepackage[utf8]{inputenc}
\usepackage{amsmath,verbatim}
\usepackage{amssymb,amsthm}
\usepackage{graphicx}
\usepackage{color}
\usepackage{enumerate,cite}
\usepackage{caption}
\usepackage{subcaption}
\usepackage{placeins}
\usepackage[margin=1in]{geometry}

\newtheorem{theorem}{Theorem}[section]
\newtheorem{lemma}[theorem]{Lemma}
\newtheorem{proposition}[theorem]{Proposition}
\newtheorem{corollary}[theorem]{Corollary}

\newenvironment{remark}[1][Remark.]{\begin{trivlist}
\item[\hskip \labelsep {\bfseries #1}]}{\end{trivlist}}

\DeclareMathOperator*{\argmin}{arg\,min}

\newcommand{\bulk}{}
\newcommand{\surf}{\textup{\text{surf}}}

\newcommand{\R}{\mathbb{R}}

\newcommand{\Z}{\mathbb{Z}}
\newcommand{\N}{\mathbb{N}}
\newcommand{\E}{\mathcal{E}}

\newcommand{\U}{\mathcal{U}}
\newcommand{\Q}{\mathcal{Q}}
\newcommand{\eps}{\varepsilon}

\newcommand{\ds}{\newline\newline\noindent}

\def\del{\delta}
\def\ddel{\delta^2}
\def\<{\langle}
\def\>{\rangle}
\def\Ea{\E^{\rm a}}

\def\Vs{V^{\surf}}
\def\Vb{V}
\newcommand{\smfrac}[2]{{\textstyle \frac{#1}{#2}}}

\title{Analysis of a Predictor-Corrector Method for Computationally Efficient
  Modeling of Surface Effects in 1D \thanks{AB acknowledges support from the DOD
    through the NDSEG Fellowship Program.  ML was supported in part by NSF PIRE
    Grant OISE-0967140, NSF Grant 1310835, and the Radcliffe Institute for
    Advanced Study at Harvard University. CO was supported by ERC Starting Grant
    335120.  }} \author{Andrew J. Binder, Mitchell Luskin, and Christoph Ortner}

\begin{document}

\maketitle

\begin{abstract}
  The regular Cauchy--Born method is a useful and efficient tool for analyzing
  bulk properties of materials in the absence of defects.  However, the method
  normally fails to capture surface effects, which are essential to determining
  material properties at small length scales.  In this paper, we present a
  corrector method that improves upon the prediction for material behavior from
  the Cauchy--Born method over a small boundary layer at the surface of a 1D
  material by capturing the missed surface effects.  We justify the separation
  of the problem into a bulk response and a localized surface correction by
establishing an error estimate, which vanishes in the long wavelength limit.
  %
\end{abstract}

\section{Introduction}

Long-range elastic fields dictate the behavior of crystalline materials in the
absence of defects.  This elastic response in crystalline systems is often
studied through the use of the Cauchy--Born method \cite{s00205-02-0218-5,
  TadmorMultiscaleBook}.  In the Cauchy--Born method, the non-local
interactions in the atomistic system are replaced by a localized continuum
approximation, which may then be further approximated using the finite element
method (FEM).  This sequence of approximations allows for a reduction in the
number of degrees of freedom in the model and yields a computationally efficient
and accurate approach to studying the material behavior of perfect crystalline
materials.  However, this approach fails in the presence of defects since they
create rapid variations in the atomistic displacement field which can no longer
be accurately captured by the Cauchy--Born model.  In particular for the present
paper, the Cauchy--Born method often struggles with capturing effects due to the
inclusion of surfaces in the model \cite{NME:NME1754}.

Surfaces can be characterized as defects because atoms near a surface
experience a different interaction environment than atoms in an idealized
crystal lattice.
For large systems, the surface influence is often negligible since the surface
region constitutes such a small portion of the entire system.  The bulk --- or
interior --- behavior dominates.  Surface effects become increasingly
significant, though, as the surface area to volume ratio increases, resulting in
a size-dependent response in material behavior, which is an active area of
research \cite{NMAT:NMAT977, PhysRevLett.95.255504, doi:10.1021/nl048456s,
  Park20083249}.


While surface effects tend to be most relevant in small systems, these small
systems may still be large enough to present computational challenges for a
fully-atomistic numerical simulation. Further, even in bulk crystals, if effects
of interest take place at or near a crystal boundary (e.g., nano-indentation),
then accurately capturing the surface physics remains crucial. Thus, a
computationally efficient means to simulate surface-dominated systems is
desirable.  Various methods have been proposed that have the potential to meet
this need.  These methods generally entail a modification of the continuum model
\cite{NME:NME1754, 2013-ac.scb.1d, BF00261375} or a concurrent coupling of
atomistic and continuum models \cite{EBTadmor:1996, 2010-MMS-SharpStabQCF,
  2010-ARMA-QCF, 2010-JMPS-qce.stab, m2an:2009007, 2011-MCOM-qnl1d, acta.atc}.
Yet, while these methods do result in a superior approximation than that which
comes from the regular Cauchy--Born method, the modified continuum models lack a
systematic control over the error in the approximation while the
atomistic-to-continuum coupling methods may encounter difficulties with surface
geometries that would result in a much-reduced computational savings than might
be observed with other defects.

In this paper, we introduce a novel approach to efficiently and accurately
capturing surface behavior through the use of a predictor-corrector method that
possesses a more controllable error and an improved approach to handling surface
geometries.  In this predictor-corrector method, the Cauchy--Born method serves
as the initial predictor for material behavior, allowing for the usual tools
that study bulk behavior to be applied to the surface problem too.  As the
Cauchy--Born method works well for the interior, it should serve well for an
approximation of the bulk response of the material.  The solution for the
Cauchy--Born method is then corrected in the next step to take into account
surface effects by minimizing the energy difference between the atomistic energy
and the Cauchy--Born energy about the Cauchy--Born solution.  This correction
occurs over a boundary layer near the surface at atomistic resolution.  The size
of the boundary layer is a controllable parameter in this predictor-corrector
method and allows for a systematic control over the accuracy of the
approximation.  The corrected solution represents the approximation of this
predictor-corrector method to the atomistic behavior.  Proving the validity of
the decomposition of the atomistic solution into a bulk response and surface
correction will be the primary goal of the analysis in this paper in addition to
assessing the quality of the approximation.

The paper is organized as follows. In \S~\ref{sec:results}, we summarize all main
results and illustrate them via numerical tests. Complete proofs are given in
\S~\ref{sec:proofs:atm}--\S~\ref{sec:proofs-pc}.

{\bf Notation for derivatives: }
We will employ three types of derivatives. (1) If $u : \{ 0, 1, \dots\}$ is a
discrete function, then we define $u_\ell' = u_{\ell+1}-u_\ell$. (2) If
$u : [0, \infty)$ is a continuous displacement or deformation function, then
$\nabla u(x)$ is the standard pointwise or weak derivative (if it exists). (3)
Finally, if $V : \R^n \to \R^p$ is a potential function (or its derivative),
then we denote its partial derivatives by $\partial_i V$, while its Jacobi
matrix is denoted by $\partial V$.


\section{Summary of Results}
\label{sec:results}

\subsection{Atomistic Model}
We consider a semi-infinite, 1D chain of atoms. We index the atoms in the chain
by the set of non-negative integers, $\Z_{\geq 0}$, and we denote the
individual location of the atom with index $\ell$ in the chain by
$y_{\ell} \in \R$.  The reference position of each atom in the chain is chosen
to be $y_{\ell} = \ell$.
We denote the displacement of atom $\ell$ from its reference position by
$u_{\ell} := y_{\ell} - \ell \in \R$.  The strain (gradient) in the bond between
the atoms indexed by $\ell$ and $\ell + 1$ in the chain will be written as
\begin{equation}\label{Eq:StrainDefinition}
 u'_{\ell} := u_{\ell + 1} - u_{\ell} = y_{\ell + 1} - y_{\ell} - 1.
\end{equation}
%
The surface of the chain is located at the atom with index 0.

The atoms in this semi-infinite chain interact according to a nearest-neighbor
site energy (effectively second-neighbour interaction). This is the simplest
setting within which we can still observe the surface effects we are interested
in. The energy due to these interactions is given by
\begin{equation}
 \E^{\rm a}(y) := V^{\surf}(y_{1} - y_{0}) +
   \sum_{\ell=1}^{\infty}\Vb(y_{\ell} - y_{\ell - 1}, y_{\ell+1} -
y_{\ell}),
\end{equation}
where $\Vb$ denotes the site energy for atoms in the interior of the
system while $V^{\surf}$ denotes the surface site energy.
%
%
The surface atom merits a different site energy than the interior atoms because
the 0-th atom has only one neighbor while every other atom has two. We will
later assume that $V^{\rm surf}$ is the limit of $\Vb$ as one of the
bonds is stretched to infinity.

Since $\E^{\rm a}$ is translation invariant, configurations are not meaningfully
distinct if they differ only up to a translation. Hence, we fix the end point so
that $y_{0} = 0$. Under this constraint, knowledge of $u'$ allows us to recover
the full displacement $y$. Hence, it will
%
%
often be more convenient to consider the energy in terms of the strain rather
than the deformation:
%
%
\begin{equation}\label{Eq:GeneralizedAtomisticEnergy}
\E^{\rm a}(u) := V^{\surf}(u_{0}') +
   \sum_{\ell=1}^{\infty}\Vb(u_{\ell - 1}', u_{\ell}'),
\end{equation}
where we have now also absorbed the reference strain $1$ into the definitions of
$V^{\surf}$ and $\Vb$.

We assume throughout that the site energies satisfy the following properties:
\ds{\bf Site Energy Properties:}
\begin{enumerate}[(i)]
\item\label{Eq:SmoothnessSiteEnergy} $\Vb \in C^{k}(\R^{2})$ and
  $V^{\surf} \in C^{k}(\R)$ with $k \geq 3$;
\item\label{Eq:EnergyAndDerivativeBounds} $\Vb, V^{\surf}$ and all
  permissible partial derivatives are bounded.
\item\label{Eq:MinSiteEnergy}
   $\Vb(0,0) = 0$;
\item\label{Eq:PositiveDefiniteSiteEnergy}
   $\partial^{2}\Vb(0,0) > 0$;
\item\label{Eq:PositiveSiteEnergy}
   $\inf_{\{|(r,s)| > \eps\}}\Vb(r,s) > 0$ for any $\eps > 0$;
\item\label{Eq:DefSurfSiteEnergy}
   For any $s \in \R$, $\lim_{r \to \infty}\Vb(r, s) = V^{\surf}(s)$;
\end{enumerate}
If the semi-infinite chain were extended infinitely in both directions so that
it had no surface and its energy consisted purely of a sum of bulk site
energies, then Properties
\eqref{Eq:MinSiteEnergy}--\eqref{Eq:PositiveSiteEnergy} would imply that the
ground state of this bulk model is the configuration with a uniform strain of 0.
Property \eqref{Eq:PositiveDefiniteSiteEnergy} would guarantee that the phonon
frequencies of such an infinite chain system are positive so that the infinite
crystal is stable.  For the surface model, these properties will allow us to
prove existence of a ground state and to establish properties of the boundary
layer.
%
Property \eqref{Eq:DefSurfSiteEnergy} defines the surface site energy in terms
of the limiting behavior of the bulk site energy, and as a consequence of
properties \eqref{Eq:PositiveSiteEnergy} and \eqref{Eq:DefSurfSiteEnergy}, we
see that
\begin{equation}\label{Eq:PositiveSurfaceSiteEnergy}
 \inf_{s \in \R} V^{\surf}(s) > 0.
\end{equation}
Property \eqref{Eq:EnergyAndDerivativeBounds} is one of convenience and is not
strictly necessary.  Throughout the paper, derivatives of $V^{\surf}$
with respect to its argument will be denoted by $\partial_{F}V^{\surf}$.
Derivatives of $V^{\bulk}$ with respect to its first and second argument will be
denoted by $\partial_{1}V^{\bulk}$ and $\partial_{2}V^{\bulk}$, respectively.
Higher-order derivatives for the two site energies will be denoted similarly.

We are concerned with finding the energy-minimizing configuration of the chain
in the presence of external forces. To consider the situation where bulk and
surface effects are roughly of the same order, we consider only {\em finite-energy
displacements}, i.e., displacements from the space
%
%
\begin{equation}\label{Def:U}
 \U := \{u : \Z_{\geq 0} \to \R \; | \; u(0) = 0 \; \text{and} \;
                                     u' \in \ell^{2}(\Z_{\geq 0})\}.
\end{equation}
When $\U$ is equipped with the $H^{1}$-seminorm
$|u|_{H^{1}} = \|u'\|_{\ell^{2}(\Z_{\geq 0})}$ it becomes a Hilbert space due to
the constraint $u(0) = 0$. It is easy to see \cite{acta.atc} that compact
displacements are dense in $\U$.

Applied forces take the form of a lattice function $f : \Z_{\geq 0} \to \R$,
which must be an element of $\U^{*}$.  We say that $f \in \U^{*}$ if and only if
there exists a constant $\|f\|_{\U^{*}}$ such that
\begin{equation}
 \langle f, u \rangle
   \leq
 \|f\|_{\U^{*}}\|u'\|_{\ell^{2}(\Z_{\geq 0})}
 \quad \text{for all} \; u \in \U \text{ with } {\rm supp}(u) \text{ compact},
\end{equation}
where $\langle f, u \rangle := \sum_{\ell = 0}^{\infty}f_{\ell}u_{\ell}$. For
arbitrary $u \in \U$, $\langle f, u \rangle$ is defined by continuity.

Given such a force $f \in \U^{*}$, we seek a minimizer
\begin{equation}\label{SemiInfiniteAtomisticProblem}
  u^{\rm a}
    \in
  \argmin
    \{ \E^{\rm a}(u) - \langle f, u\rangle_{\Z_{\geq 0}} \; | \; u \in \U \}.
\end{equation}
This problem may have several or no solutions.  We consider the existence of
solutions to this problem as well as their decay in the following section.

We conclude the description of the model by mentioning that Properties
\eqref{Eq:SmoothnessSiteEnergy} and \eqref{Eq:EnergyAndDerivativeBounds}
guarantee that the first and second variations are globally Lipschitz
continuous.  We refer to \cite{acta.atc} for a proof.


\begin{remark}
  It is easy to see that $f \in \U^*$ if and only if there exists
  $g \in \ell^2(\Z_{\geq 0})$ such that $f = -g'$ by considering a discrete
  summation by parts of $\langle f, u \rangle$.  From this summation by parts,
  it can also be shown that
  \begin{equation}
   g_{n} = \sum_{\ell = n + 1}^{\infty}f_{\ell} < \infty;
  \end{equation}
  see \S~\ref{sec:proofs-cb} for a proof that this sum is well-defined.
  Therefore,
  $\langle f, u \rangle = \langle g, u' \rangle \leq \| g \|_{\ell^2} \| u'
  \|_{\ell^2}$, and in fact, it is clear that $\|f\|_{\U^{*}} = \|g\|_{\ell^{2}}$.


\end{remark}

\subsubsection{Example: EAM Model}
To provide a concrete example of the types of systems encompassed by the model
described above, we show how the semi-infinite chain satisfying the Embedded
Atom Model (EAM) \cite{PhysRevB.37.3924,PhysRevB.29.6443} with only
nearest-neighbor interactions is described in this framework.  The energy for
such a system is written in terms of the strain as
\begin{align}\label{Eq:EAMModel}
 \E^{\rm a}(u') &= \phi(1 + u'_{0}) + \psi(\rho(1 + u'_{0}))
+ \sum_{\ell=1}^{\infty}\biggl[
        \phi(1 + u'_{\ell})
      + \psi\Big(\rho(1 + u'_{\ell}) + \rho(1 + u'_{\ell-1})\Big)
    \biggr],
\end{align}
where $\phi$ is a nearest-neighbor pair potential, $\psi$ is an embedding energy
function, and $\rho$ is an electron density function.
In order to write this energy in the form of
\eqref{Eq:GeneralizedAtomisticEnergy}, we define the bulk and surface site
energies, respectively, as
\begin{align*}
 \Vb(u_{\ell - 1}', u_{\ell}') &=
   \frac{1}{2}\phi(1 + u_{\ell - 1}') + \frac{1}{2}\phi(1 + u_{\ell}') +
   \psi\bigl(\rho(1 + u_{\ell -1}') + \rho(1 + u_{\ell}')\bigr),
                                       \quad \text{and} \\
 V^{\surf}(u_{0}') &=
   \frac{1}{2}\phi(1 + u_{0}') + \psi\bigl(\rho(1 + u_{0}')\bigr).
\end{align*}
For the analysis, we will consider the generalized framework described by
\eqref{Eq:GeneralizedAtomisticEnergy}, but we will return to the EAM system for
the numerical results in \S~\ref{sec:numerical-results}.

\subsection{Existence, Decay and Stability of Atomistic Solutions}
Unlike in a bulk model without surfaces, we do not expect that an
energy-minimizing configuration for the semi-infinite chain is
homogeneous. Therefore, we  have to be satisfied with weaker
existence proofs and establishing general properties of a minimizer.

\begin{theorem}\label{Thm:ExistenceMinimizerAtomistic}
  There exists a minimizer of $\E^{\rm a} : \U \to \R \cup \{+\infty\}$.
\end{theorem}
\begin{proof}
  See Section \ref{sec:proofs:atm}.
\end{proof}
While for many systems it is reasonable (and natural) to expect a unique ground
state, our assumptions do not preclude the existence of multiple states that
achieve the same minimal energy. In the following, $u^{\rm a}_{\rm gr}$ may
refer to any ground state.

The key property of $u^{\rm a}_{\rm gr}$ that motivates our subsequent
developments is that the surface effects are highly localized. This is
established next.
\begin{theorem}\label{Thm:ExpDecayAtomisticStrain}
  Let $u^{\rm a}_{\rm gr}$ be a critical point of $\E^{\rm a}$.  Then, there
  exists $0 \leq \mu_{\rm a} < 1$ such that
  \begin{equation}
    |(u^{\rm a}_{\rm gr})'_{\ell}|
      \lesssim
    \mu_{\rm a}^{\ell}
    \qquad \text{for all} \quad \ell \in  \Z_{\geq 0}.
  \end{equation}
\end{theorem}
\begin{proof}
  See Section \ref{sec:proofs:atm}.
\end{proof}

The exponential decay of the strain due to surface effects normally justifies
ignoring surface effects in large systems. However, we will see that when
surface effects are the focus of interest, then this boundary layer cannot be
ignored.



We now proceed to incorporating external forces into the analysis.  To that end,
it is convenient to assume that a ground state $u^{\rm a}_{\rm gr}$ is strongly
stable; that is, we suppose that there exists an {\em atomistic stability
constant} $c_{\rm a} > 0$ such that
\begin{equation}\label{Eq:AtomisticGroundStateStability}
 \langle \delta^{2}\E^{\rm a}(u^{\rm a}_{\rm gr})v, v \rangle
   \geq
 c_{\rm a}\|v'\|^{2}_{\ell^{2}(\Z_{\geq 0})}
   \quad \text{for all} \; v \in \U.
\end{equation}
This stability assumption on $u^{\rm a}_{\rm gr}$ enables us to prove the
existence of nearby strongly stable local minimizers of the atomistic problem
from \eqref{SemiInfiniteAtomisticProblem} with small external forces. For future
reference, an element $u^{\rm a} \in \U$ is a strongly stable solution to the
atomistic problem if and only if it satisfies the Euler--Lagrange equation
\begin{equation}
 \langle \delta \E^{\rm a}(u^{\rm a}), v \rangle
 =
 \langle f, v \rangle \quad \text{for all} \; v \in \U
\end{equation}
as well as the stability condition
\begin{equation}\label{Eq:AtomisticStability}
 \langle \delta^{2} \E^{\rm a}(u^{\rm a})v, v\rangle
   \geq
 c\|v'\|^{2}_{\ell^{2}(\Z_{\geq 0})}
   \quad \text{for all} \; v \in \U,
\end{equation}
for some constant $c > 0$.  The exact form of the first and second variations of
$\E^{\rm a}$ are provided in Propositions \ref{Prop:FirstVariations} and
\ref{Prop:SecondVariations}.
\begin{corollary}\label{Thm:AtomisticLocalMinimizers}
  There exist $\eps, C > 0$ such that, for all $f \in \U^{*}$ with
  $\|f\|_{\U^{*}} < \eps$, the atomistic problem
  \eqref{SemiInfiniteAtomisticProblem} has a unique, strongly-stable solution
  with $\| (u^{\rm a} - u^{\rm a}_{\rm gr})' \|_{\ell^2} \leq C \| f \|_{\U^*}$.
\end{corollary}
\begin{proof}
  This is an immediate consequence of the inverse function theorem (Theorem
  \ref{Thm:InverseFunctionTheorem}).
\end{proof}
%

\subsection{Cauchy--Born Model}

A common approach to determining the approximate bulk behavior of perfect
crystalline materials utilizes the Cauchy--Born model of atomistic interactions,
approximates the non-local atomistic interactions with a local approximation.  A
limiting procedure then turns the discrete problem into a continuum one
\cite{s00205-02-0218-5}.  The Cauchy--Born energy for the semi-infinite chain
model is
\begin{equation}\label{Eq:CauchyBornEnergy}
 \E^{\rm cb}(u)
   :=
 \int_{0}^{\infty}W(\nabla u(x)) \,dx
   \qquad \text{for}\; u \in \U^{\rm cb},
\end{equation}
where $W(F) := V^{\bulk}(F, F)$ is the Cauchy--Born energy density function and
\begin{equation}
  \U^{\rm cb}
    :=
  \big\{u \in H^{1}_{\text{loc}}(0,\infty) \; | \; \nabla u \in L^{2}(0,\infty)
    \; \text{and} \; u(0) = 0 \big\}.
\end{equation}
The Cauchy--Born energy density function $W$ inherits the smoothness of
$V^{\rm bulk}$, i.e., $W \in C^{3}(\R)$.
Derivatives of $W$ with respect to its argument (as opposed to $x$) will be
indicated by $\partial_{F}W$.  We also note that the space $\U$ may be
considered a subspace of $\U^{\rm cb}$ if we identify the lattice functions
$u \in \U$ with their piecewise continuous interpolants.

For atomistic systems without defects, the Cauchy--Born approximation and the
exact model agree under homogeneous deformation.  However, in systems containing
defects, such as surfaces, this property is lost.  In
\eqref{Eq:CauchyBornEnergy}, this discrepancy arises as the Cauchy--Born model
treats every point in the chain as an interior point.  The absence of a surface
component to the energy results in an $\mathcal{O}(1)$ error in the consistency
estimate for the Cauchy--Born model as compared to the atomistic system, which
we can demonstrate when $f = 0$.

%
%
\begin{proposition}\label{Thm:CBConsistencyErrorAtSurface}
  The unique minimizer of $\E^{\rm cb}$ in $\U^{\rm cb}$ is $u^{\rm cb} =
  0$. Its atomistic residual is bounded by
  \begin{equation}
    \sup_{v \in \U, \|v'\|_{\ell^{2}(\Z_{\geq 0})} = 1}
    |\langle \delta \E^{\rm a}(0), v \rangle |
    =
    \|\delta\E^{\rm a}(0)\|_{\U^{*}}
    =
    |\partial_{F}V^{\surf}(0)|.
  \end{equation}
  In particular,
  $\| (u^{\rm a}_{\rm gr} - u^{\rm cb})' \|_{\ell^2} \geq M^{-1}
  |\partial_{F}V^{\surf}(0)|$,
  where $M$ is the global Lipschitz constant of $\delta \E^{\rm a}$.
\end{proposition}
\begin{proof}
  See Section \ref{sec:proofs-cb}.
\end{proof}
In general, we expect that $|\partial_{F}V^{\surf}(0)| \neq 0$ so that the
Cauchy--Born and atomistic energy-minimizing configurations are not in agreement
because of surface effects.


This lack of consistency at the surface indicates that the Cauchy--Born model
alone will not serve as an accurate approximation of the original system for the
purpose of studying surface effects.  However, as we have shown in Theorem
\ref{Thm:ExpDecayAtomisticStrain}, surface effects are a local phenomenon.  The
Cauchy--Born model may then still serve as an efficient means to computing the
behavior of the system in the interior.

We will show that the bulk and surface effects decouple, which means that a
complex concurrent coupling scheme is not required. Instead, we can ``add'' the
surface effects in a {\em predictor-corrector type approach}, where the
predictor is the Cauchy--Born solution.


We therefore consider the general Cauchy--Born problem
\begin{equation}\label{SemiInfiniteCauchyBornProblem}
 u^{\rm cb} \in \argmin\{\E^{\rm cb}(u) - \langle f, u \rangle_{\R_{+}} \; | \;
                         u \in \U^{\rm cb}\},
\end{equation}
where we identify the lattice function $f$ with its continuous piecewise affine
interpolant and $\langle f, u \rangle_{\R_{+}} := \int_{0}^{\infty} f u\, dx$.
It is easy to see that $f \in \U^*$ implies that
$\langle f, \cdot \rangle_{\R_+} \in (\U^{\rm cb})^*$; see
\S~\ref{sec:proofs-cb}.


Analogously to the atomistic problem, we say $u^{\rm cb}$ is a strongly stable
solution to \eqref{SemiInfiniteCauchyBornProblem} if it satisfies the
first-order and strong second-order optimality conditions:
\begin{align}
  \notag
  \langle \delta \E^{\rm cb}(u^{\rm cb}), v \rangle
  &=
    \langle f, v \rangle_{\R_{+}} \quad \text{for all} \; v \in \U^{\rm cb}
    \qquad \text{and} \\
  \label{Eq:CBStability}
  \langle \delta^{2} \E^{\rm cb}(u^{\rm cb})v, v\rangle
  &\geq c_{\rm cb}\|\nabla v\|^{2}_{L^{2}([0,\infty))}
    \quad \text{for all} \; v \in \U^{\rm cb}
\end{align}
for some constant $c_{\rm cb} > 0$.  The exact forms of the first and second
variations of the Cauchy--Born energy may be found in Propositions
\ref{Prop:FirstVariations} and \ref{Prop:SecondVariations}.

For small enough external forces, we can guarantee the existence of strongly
stable solutions to the Cauchy--Born problem and deduce some additional facts
concerning their regularity.
\begin{theorem}\label{Thm:CBExistenceRegularity}
  There exists $\eps_{\rm cb} > 0$ such that for all $f \in \U^{*}$ with
  $\|f\|_{\U^{*}} < \eps_{\rm cb}$, a strongly stable solution
  $u^{\rm cb} \in \U^{\rm cb}$ of the Cauchy--Born problem
  \eqref{SemiInfiniteCauchyBornProblem} exists.

  Moreover, $u^{\rm cb} \in H^3_{\rm loc}$, and it satisfies the bounds
  \begin{equation}
    |\nabla u^{\rm cb}(0)| \lesssim \|f \|_{\U^*}, \quad
    |\nabla^{2}u^{\rm cb}(x)| \lesssim |f(x)|
    \quad
    \text{and}
    \quad
    |\nabla^{3}u^{\rm cb}(x)| \lesssim |\nabla f(x)| + |f(x)|^{2}.
  \end{equation}
  %
\end{theorem}
\begin{proof}
  See Section \ref{sec:proofs-cb}.
\end{proof}
%


In the next section, we will use the Cauchy--Born solution to approximate the
atomistic solution. To that end, we introduce a projection operator
$\Pi_{\rm a} : \U^{\rm cb} \to \U$ via
\begin{displaymath}\label{Eq:CBProjector}
  (\Pi_{\rm a} u)_0 = 0 \qquad \text{and} \qquad
  (\Pi_{\rm a} u)'_\ell = \int_{\ell}^{\ell+1} \nabla u(s) \,ds
  \quad \text{for } \ell \in \Z_{\geq 0}, u \in \U^{\rm cb}.
\end{displaymath}

\subsection{Predictor-Corrector Method}
Proposition \ref{Thm:CBConsistencyErrorAtSurface} demonstrates that the
Cauchy--Born method commits an $O(1)$ consistency error at a material
surface. However, the Cauchy--Born method is an excellent method for
approximating the bulk response of materials, and Theorem
\ref{Thm:ExpDecayAtomisticStrain} indicates that the error due to surface
effects may be quite localized.  We therefore propose a mechanism to {\em
  correct} the Cauchy--Born solution to obtain an approximation of the form
\begin{equation} \label{Eq:ApproxDecompOfAtomisticStrain}
  u^{\rm a} \approx \Pi_{\rm a} u^{\rm cb} + q_L,
\end{equation}
where the parameter $L$ for the corrector $q_{L}$ determines the quality
of the approximation.

Given a {\em predictor} $u^{\rm cb}$ solving the Cauchy--Born problem
\eqref{SemiInfiniteCauchyBornProblem}, we define the corrector problem via
the minimization of a {\em corrector energy}, which is given by
\begin{align}\label{Eq:CorrectorEnergy}
 \E^{\Gamma}(q; F_{0}) &=
   V^{\surf}(F_{0} + q_{0}')
     - W(F_{0}) - q_{0}'\partial_{F}W(F_{0})
 \\ & \qquad +
      \notag
  \sum_{j=1}^{\infty}\left(
    \Vb(F_{0} + q_{j-1}', F_{0} + q_{j}')
    - W(F_{0})
    - q_{j}'\partial_{F}W(F_{0})\right),
\end{align}
where we will normally take $F_{0} := \nabla u^{\rm cb}(0)$.

The idea of the corrector problem is that it should depend only in a {\em local}
way on the elastic field $\nabla u^{\rm cb}$.  In this case, this dependence is
only on $\nabla u^{\rm cb}(0)$. This choice was made deliberately so that the corrector
problem can be understood as a cell problem on a surface element when the
Cauchy--Born model is discretized using finite elements.


For $L \in \N \cup \{\infty\}$, a corrector strain on the interval $[0, L]$ is
found by solving
%
\begin{equation}\label{CorrectorProblem}
  q_{L} \in \argmin\{\E^{\Gamma}(q; F_{0}) \; | \; q \in \Q_{L}\},
\end{equation}
where $\Q_{L}$ describes the the boundary layer over which we correct the
Cauchy--Born solution and is defined to be
\begin{align*}
  \Q_{L} &:= \{q \in \U \; | \; q_{\ell}' = 0 \; \text{for all} \; \ell \geq L\}.
           \quad \text{In particular,}  \quad Q_\infty = \U.
\end{align*}
The choice of $L$ affects the computational expense to solve the corrector
problem as well as the accuracy of the approximation of the atomistic system's
behavior.  Note that the external force $f$ does not enter directly into the
corrector problem.  The Cauchy--Born method accounts for the external forces on
its own, but the influence of the external force is felt through the dependence
on $F_{0}$.

%
\begin{theorem}\label{Thm:ExistenceInfiniteCorrector}
  There exists $\eps_{\Gamma} > 0$ such that, for all
    $F_0 \in \R$ with $|F_0| < \eps_\Gamma$, the corrector
  problem~\eqref{CorrectorProblem} with $L = \infty$ has a solution
  $q_{\infty} \in \Q_{\infty}$.  For all $v \in \Q_{\infty}$, this solution
  $q_{\infty}$ satisfies
\begin{equation}
  \label{eq:corrector-stab}
 \langle \delta^{2} \E^{\Gamma}(q_{\infty}; F_{0})v, v \rangle
   \geq \frac{c_{\rm a}}{2}\|v'\|_{\ell^{2}}^{2},
\end{equation}
where $c_{\rm a}$ is the atomistic stability constant for $u^{\rm a}_{\rm gr}$ in
\eqref{Eq:AtomisticGroundStateStability}.  In addition, there exists a constant
$0 \leq \mu_{q} < 1$ such that
\begin{equation}
  \label{eq:corrector-exp-decay}
  |(q_{\infty})'_{\ell}| \lesssim \mu_{q}^{\ell} \qquad \text{for all} \quad
    \ell \in \Z_{\geq 0}.
\end{equation}
\end{theorem}
\begin{proof}
 See Section \ref{sec:proofs-corrector}.
\end{proof}
%
We can now consider the existence of a solution to the corrector problem for a
finite boundary layer as an approximation of the infinite case.
\begin{proposition}\label{Prop:ExistenceFiniteCorrector}
  Under the conditions of Theorem \ref{Thm:ExistenceInfiniteCorrector},
  there exists $L_{0} > 0$ such that the corrector
  problem~\eqref{CorrectorProblem} with $L \geq L_0$ has a solution $q_{L} \in
  \Q_{L}$ satisfying
  \begin{displaymath}
    \|q_{\infty}' - q_{L}'\|_{\ell^{2}} \lesssim \mu_{q}^{L}.
  \end{displaymath}
\end{proposition}
\begin{proof}
 See Section \ref{sec:proofs-corrector}.
\end{proof}

  If $\|f\|_{\U^*}$ is sufficiently small, then Theorem
  \ref{Thm:CBExistenceRegularity} guarantees the existence of a solution
  $u^{\rm cb} \in \U^{\rm cb}$ with $|\nabla u^{\rm cb}(0)| \leq \eps_\Gamma$.
  Thus, if we choose a sufficiently large boundary layer for the corrector
  problem \eqref{CorrectorProblem}, Theorem \ref{Thm:ExistenceInfiniteCorrector}
  and Proposition \ref{Prop:ExistenceFiniteCorrector} ensure that the {\em
  predictor-corrector approximation}
  \begin{displaymath}
    u^{\rm pc}_L := \Pi_{\rm a}u^{\rm cb} + q_{L}
  \end{displaymath}
  is well-defined.  Next, to determine the accuracy of this approximation, we
  show that a solution to the atomistic problem
  \eqref{SemiInfiniteAtomisticProblem} exists in a neighborhood of
  $u^{\rm pc}_L$, which we quantify.
\begin{theorem}
  There exists an $\eps > 0$ such that, for all $f \in \U^{*}$ with
  $\|f\|_{\U^{*}} < \eps$, there exists an atomistic solution $u^{\rm a} \in \U$
  to \eqref{SemiInfiniteAtomisticProblem} satisfying
  \begin{equation}\label{Eq:PCApproximationError}
    \big\| (u^{\rm a})' -(u^{\rm pc}_L)' \big\|_{\ell^{2}}
    \lesssim
    \mu_{q}^{L}
    + |\nabla^2 u^{\rm cb}(0)|
    + \|\nabla^{2} u^{\rm cb}\|_{L^{4}}^{2} + \|\nabla^{3} u^{\rm cb}\|_{L^{2}}
    + \| \nabla f \|_{L^2}.
  \end{equation}
\end{theorem}
\begin{proof}
  This result is a consequence of the inverse function theorem (Theorem
  \ref{Thm:InverseFunctionTheorem}) and Theorems~\ref{Thm:PCConsistency}
  (Consistency) and \ref{Thm:PCStability} (Stability).
\end{proof}
\begin{remark}
The term $|\nabla^2 u^{\rm cb}(0)|$ is an error due to the fact that, in the
corrector problem, we replaced the nonlinear elastic field
$\nabla u^{\rm cb}(x)$ with a homogeneous field $\nabla u^{\rm cb}(0)$.
Identifying this term is the main result of our analysis.

The term $\mu_q^L$, which arises due to the approximation of the corrector
problem, may now be balanced against the (uncontrollable) $|\nabla^2 u^{\rm
 cb}(0)|$ contribution by choosing $L \approx \log |\nabla^2 u^{\rm
 cb}(0)|$. The remaining terms are the typical error committed by the
Cauchy--Born model and the error committed in the approximation of the external
forces, which are represented by $\|\nabla^{2} u^{\rm cb}\|_{L^{4}}^{2} +
\|\nabla^{3} u^{\rm cb}\|_{L^{2}}$ and $\|\nabla f\|_{L^2}$, respectively.

\end{remark}

In preparation for the next section, we observe that the error estimate can be
rewritten, using Theorem \ref{Thm:CBExistenceRegularity}, in terms of $f$ only:
\begin{equation}\label{Eq:PCErrorInExtForce}
  \big\| (u^{\rm a})' -(u^{\rm pc}_L)' \big\|_{\ell^{2}}
  \lesssim
  \mu_{q}^{L} + |f(0)| + \|f\|_{L^{4}}^{2} + \|\nabla f\|_{L^{2}}.
\end{equation}

\subsection{Rescaling of the External Forces}
\label{sec:rescaling}
We now consider a standard scaling of the external force typically employed
in the analysis of the Cauchy--Born model \cite{acta.atc, 2012-ARMA-cb}. For
the sake of simplicity, we discuss this only formally.

Let $\lambda^{-1}$ denote a length-scale over which we expect elastic strains to
vary. This suggests that after a rescaling of space through $x \leadsto \lambda x$
and $u \leadsto \lambda u$, the variation will be on an $O(1)$ scale. The
corresponding dual scaling for the force is $f(x) \leadsto \lambda f(\lambda x)$,
which motivates us to consider an external force $f \equiv f^{(\lambda)}$ given by
\begin{displaymath}
  f^{(\lambda)}_\ell :=  \lambda \hat{f}(\lambda \ell),
\end{displaymath}
where $\hat{f} = -\nabla \hat{g}$ for some $\hat{g} \in H^2(0, \infty)$ with
$\|\hat{g} \|_{H^2}$ sufficiently small.  We then obtain that
\begin{align*}
  |f^{(\lambda)}(0)| = \lambda |\hat{f}(0)|, \quad
  \| f^{(\lambda)} \|_{L^4}^2
    = \lambda^{3/2} \| \hat{f} \|_{L^4}^2,
  \quad \text{and} \quad
  \| \nabla f^{(\lambda)} \|_{L^2}
    = \lambda^{3/2} \| \nabla \hat{f} \|_{L^2}.
\end{align*}

  We remark that, in the following, we do not use rescaled $x$ and $u$, but
  we only mentioned this rescaling to motivate our chosen scaling of the
  external force.

For $\|\hat{g}\|_{H^2}$ sufficiently small and $L$ sufficiently large, we can
deduce the existence of the predictor $u^{\rm cb}$, the corrector $q_L$, and the
corresponding atomistic solution $u^{\rm a}$.  We can then arrive at the resulting
error estimate
\begin{displaymath}
  \big\| (u^{\rm a})' - (u^{\rm pc}_L)' \big\|_{\ell^{2}}
  \lesssim  \mu_q^L + \lambda + \lambda^{3/2} ,
\end{displaymath}
where we kept the $O(\lambda^{3/2})$ term to emphasize that the Cauchy--Born
contribution to the error is negligible compared to the surface contribution.
 We may now balance the error again by choosing $L = \log \lambda / \log \mu_q$
to finally obtain
\begin{displaymath}
  \big\| (u^{\rm a})' - (u^{\rm pc}_L)' \big\|_{\ell^{2}}
  \lesssim \lambda + \lambda^{3/2}.
\end{displaymath}

\begin{remark}
  If we take $f_\lambda$ as a smooth function rather than a piecewise affine
  interpolant, then the Cauchy--Born solution is rescaled due to the force
  rescaling according to
  $\nabla u^{cb}_{\lambda}(x) = \nabla u^{cb}(\lambda x)$.  In particular,
  $\nabla u^{\rm cb}_\lambda(0)$ would then be independent of $\lambda$, and
  thus the corrector problem is independent of $\lambda$ as well.

  In the setting of our analysis, our choice of interpolating $f$, means this
  observation is still approximately true. In particular, it is fairly
  straightforward to extend our analysis rigorously to the setting of
  \S~\ref{sec:rescaling}.
\end{remark}


\subsection{Numerical Results}
\label{sec:numerical-results}
%
In this section, we provide numerical demonstrations of the predictor-corrector
method as well as corroboration of the error estimates we established.

As atomistic model use the EAM site energies described in \eqref{Eq:EAMModel}
with $\phi, \psi,$ and $\rho$ given by
\begin{align*}
  \phi(r) &= \phi_{e}\text{exp}
    \biggl(-\gamma \biggl(\frac{r}{r_{e}} - 1\biggr)\biggr),
  \quad
  \rho(r) = f_{e}\text{exp}
    \biggl(-\beta \biggl(\frac{r}{r_{e}} - 1\biggr)\biggr),
  \\
  \psi(\rho) &= -E_{c}\biggl[1 - \frac{\alpha}{\beta}
    \log{\biggl(\frac{\rho}{\rho_{e}}\biggr)}\biggr]
    \biggl(\frac{\rho}{\rho_{e}}\biggr)^{\alpha/\beta} -
    \phi_{e}\biggl(\frac{\rho}{\rho_{e}}\biggr)^{\gamma/\beta},
\end{align*}
where $\phi_{e} = 10.6$, $f_{e} = 1$, $E_{c} = 3.54$, $\alpha = 21$,
$\beta = 6$, $\rho_{e} = 2$, $r_{e} = 1$, and $\gamma = 8$.  These parameters
and potentials are taken from \cite{PhysRevB.37.3924} and represent an EAM
potential describing a system composed of copper atoms.  We note that $r_{e}$
represents the equilibrium distance in the infinite atomistic model without
surfaces. Therefore, the equilibrium spacing in the infinite model is simply 1
and corresponds to a strain of 0 as our reference spacing is also 1.



For the numerical implementation, we seek energy-minimizing configurations of
the semi-infinite atomistic chain model in the space $\U_N := \Q_N$,
%
%
instead of $\U$. Provided $N$ is sufficiently large so that
${\rm supp}(f) \subset [0, N/2]$, say, it can be proven using the techniques we
employed in our analysis that the additional error committed is exponentially
small in $N$. Due to the low computational cost in one-dimensional experiments,
we can choose $N = 1000$ throughout, which guarantees that the additional error
committed from replacing $\U$ with $\U_N$ is negligible.

We find approximate solutions of the Cauchy--Born
problem~\eqref{SemiInfiniteCauchyBornProblem} by seeking minimizers of a
discretized Cauchy--Born energy given by
\begin{displaymath}
  \E^{\rm cb}(u) := \sum_{\ell = 0}^{\infty}W(u_{\ell}'),
  \quad \text{where} \quad
  W(F) =
  \frac{1}{2}\phi(1 + F) + \frac{1}{2}\phi(1 + F) +
  \psi\bigl(2\rho(1 + F)\bigr).
\end{displaymath}
%
%
Proposition \ref{Proof:CBDiscretization} can be employed to prove that this
discretization does not introduce new terms to the error bound
\eqref{Eq:PCApproximationError}.

The energy-minimizing configurations for the above models are found using the
steepest-descent method with a backtracking algorithm.

\begin{figure}
   \centering
   \begin{subfigure}[b]{0.48\textwidth}
     \includegraphics[scale=0.40]{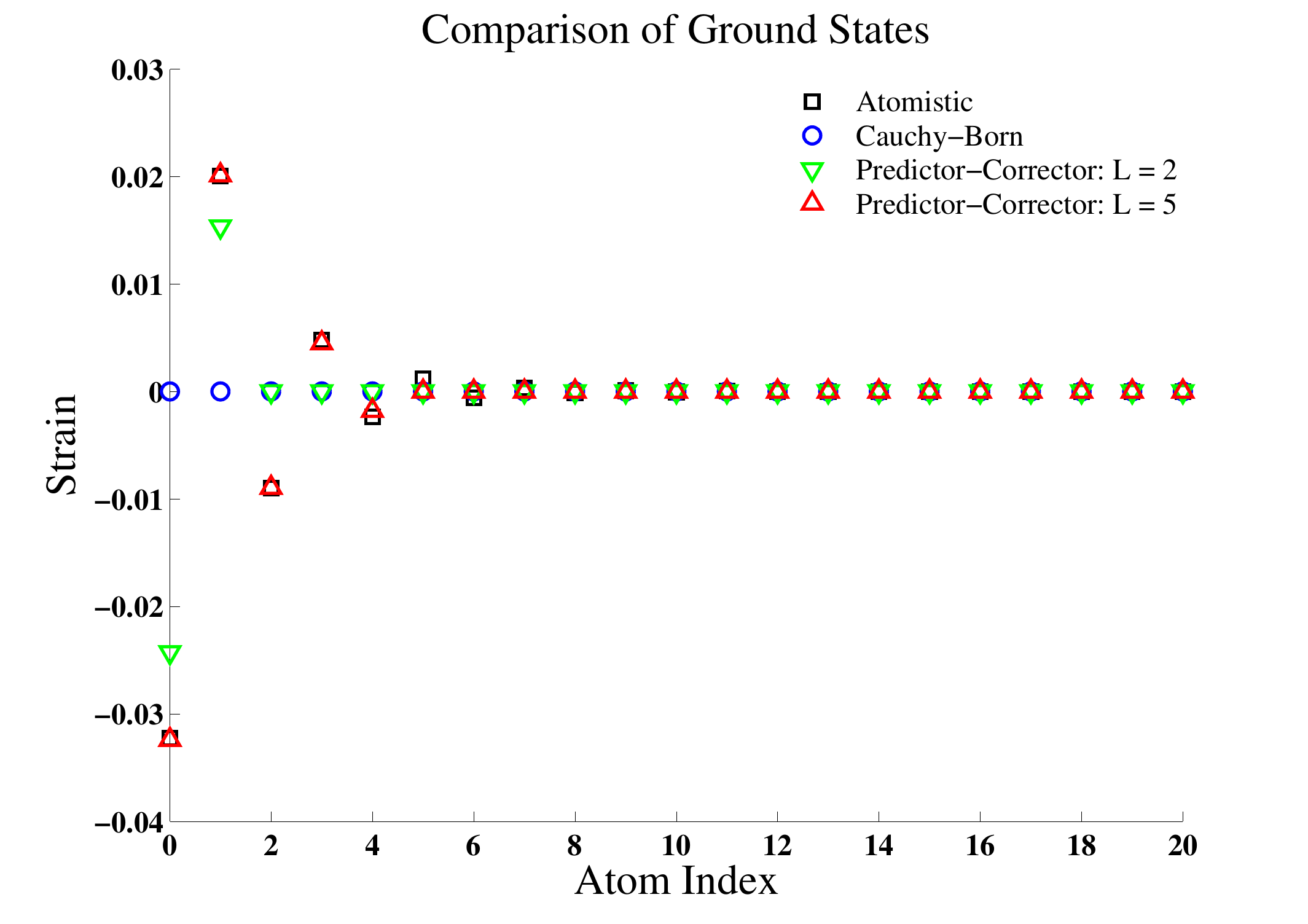}
     \caption{ }
     \label{fig:GroundStates}
   \end{subfigure}
   ~
   \begin{subfigure}[b]{0.48\textwidth}
     \includegraphics[scale=0.40]{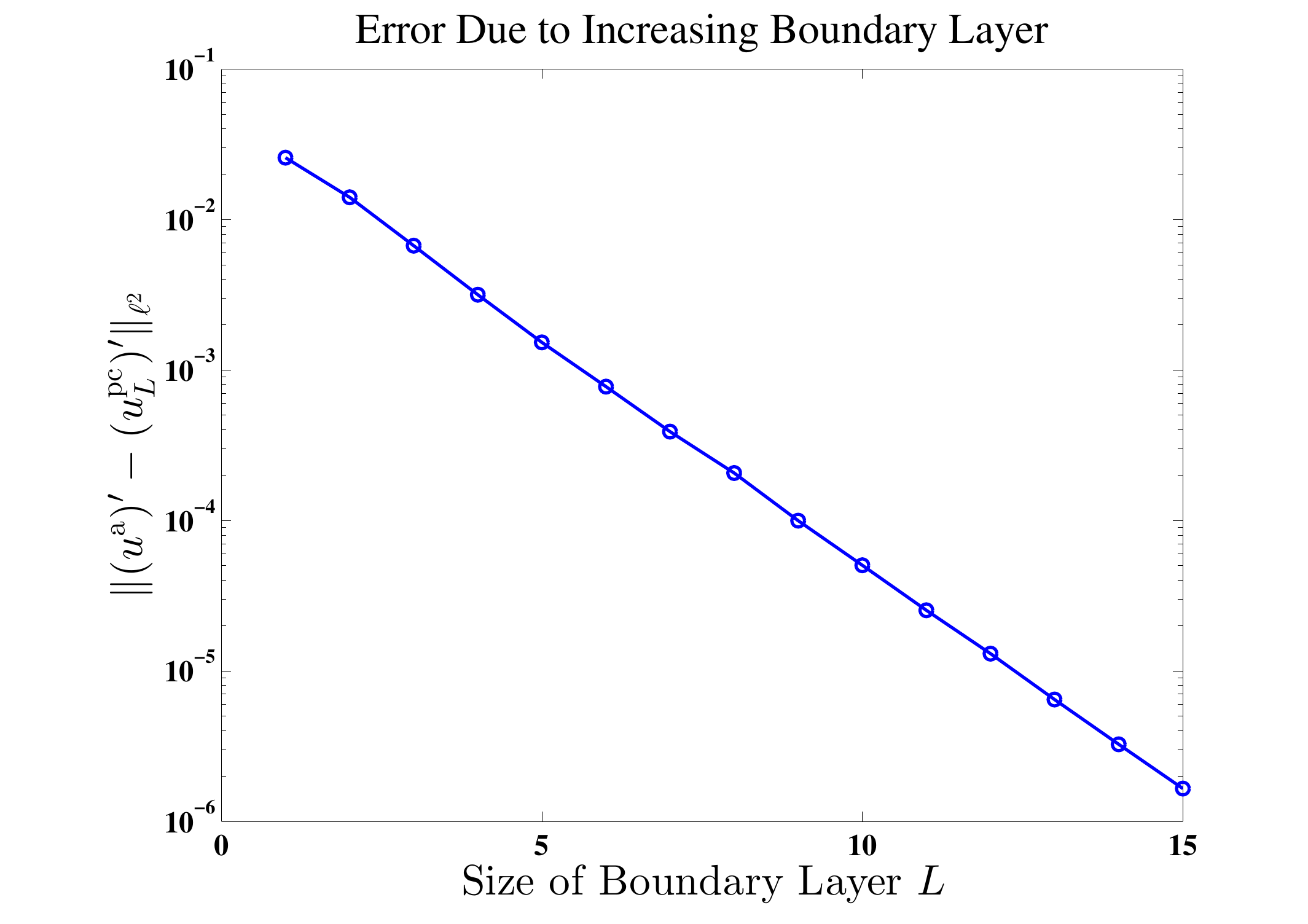}
     \caption{ }
     \label{fig:ExponentialConvergence}
   \end{subfigure}
   \caption{Ground states for the atomistic, Cauchy--Born, and
            predictor-corrector methods along with a demonstration of the
            exponential convergence of the boundary layer error.}
\end{figure}

\subsubsection{Test 1: Ground State}
The ground states ($f = 0$) of the atomistic, pure Cauchy--Born, and
predictor-corrector method are shown in Figure \ref{fig:GroundStates}.  The
atomistic ground state exhibits the predicted exponential decay, but we also
observe an alternating sign behavior common to certain metals such as copper.

The pure Cauchy--Born solution clearly does not capture any surface effects.
For the predictor-corrector method, we see that $u^{\rm pc}_2$ already yields
good accuracy while $u^{\rm pc}_5$ is visually indistinguishable from
$u^{\rm a}_{\rm gr}$.  In Figure \ref{fig:ExponentialConvergence}, we quantify
this by numerically demonstrating the exponential convergence of $q_L = u^{\rm
  pc}_L$ to the ground state.


\FloatBarrier

\subsubsection{Test 2: Long Wavelength Limit}
\label{sec:test2}
\begin{figure}
   \centering
   \includegraphics[scale=0.40]{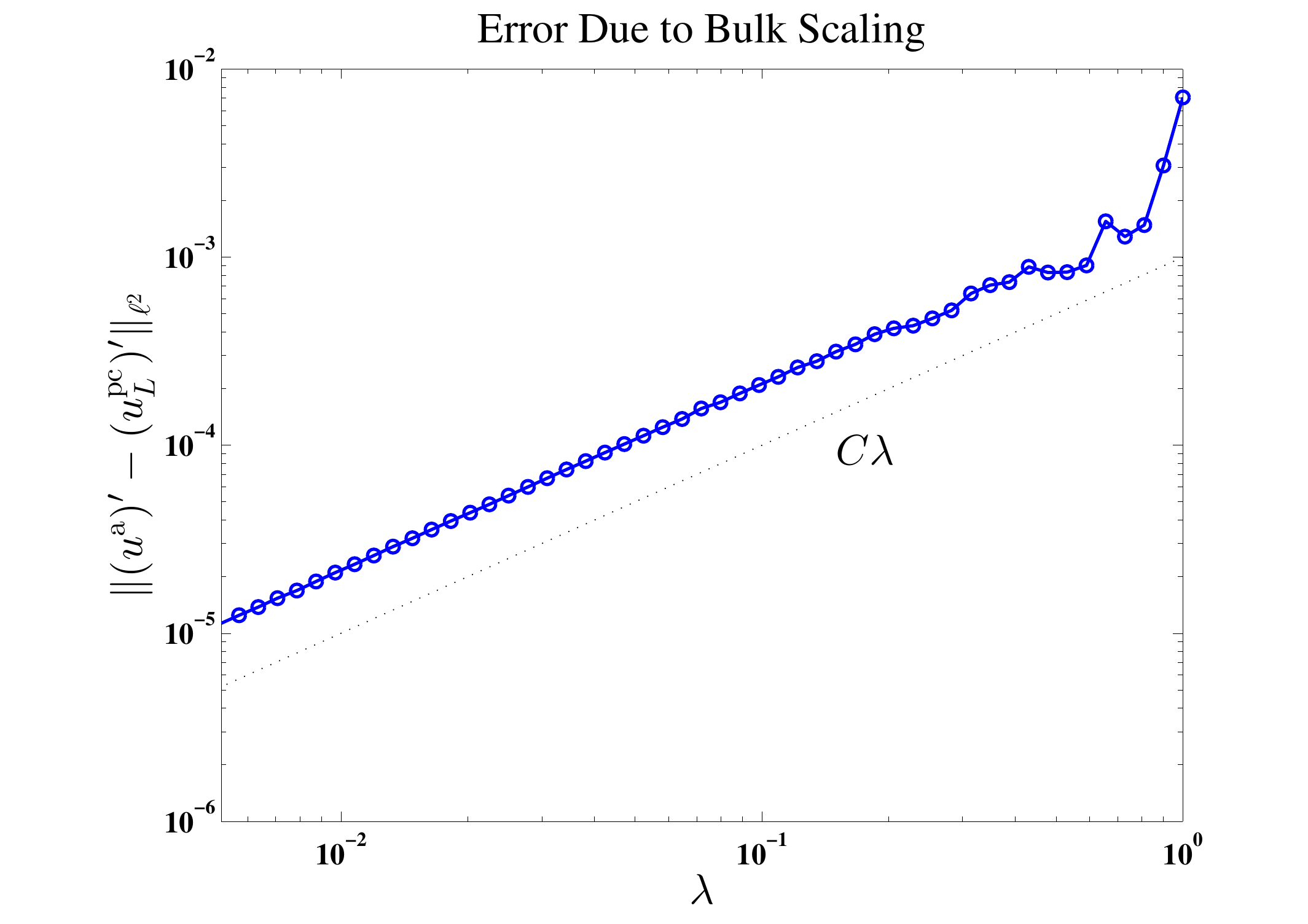}
   \caption{Rate of convergence to atomistic solution in the long-wavelength
     limit.}
    \label{fig:FirstOrderConvergence}
\end{figure}

We next consider the error under the long wavelength limit and rescaling
described in \S~\ref{sec:rescaling}.  We define
$\hat{f}(x) = \cos(3\pi x/8)\chi_{[0,4)}(x)$, where $\chi_{A}(x)$ denotes the
characteristic function, and $f_\ell = \lambda \hat{f}(\lambda \ell)$. For each
$\lambda$, we compute the solution to the predictor-corrector method with the
size of the boundary layer set to $L = 3 + \lceil \log_{10/9}(\lambda^{-1})
\rceil$.  This choice of $L$ is motivated from the
discussion of balancing the boundary layer contribution
to the error in \S~\ref{sec:rescaling}.  By choosing a logarithm with a base
less than $\mu_{q}^{-1}$ in the choice of $L$, the surface contribution to the
error becomes $\lambda^{\eta}$ where $\eta > 1$ is a constant.  Figure
\ref{fig:FirstOrderConvergence} shows the first-order convergence that results
when the scaling factor $\lambda$ is taken to zero.  The predictor-corrector
method behaves precisely as predicted in \S~\ref{sec:rescaling}.

\subsubsection{Test 3: Error with External Forces}
\begin{figure}
   \centering
   \begin{subfigure}[b]{0.48\textwidth}
     \includegraphics[scale=0.40]{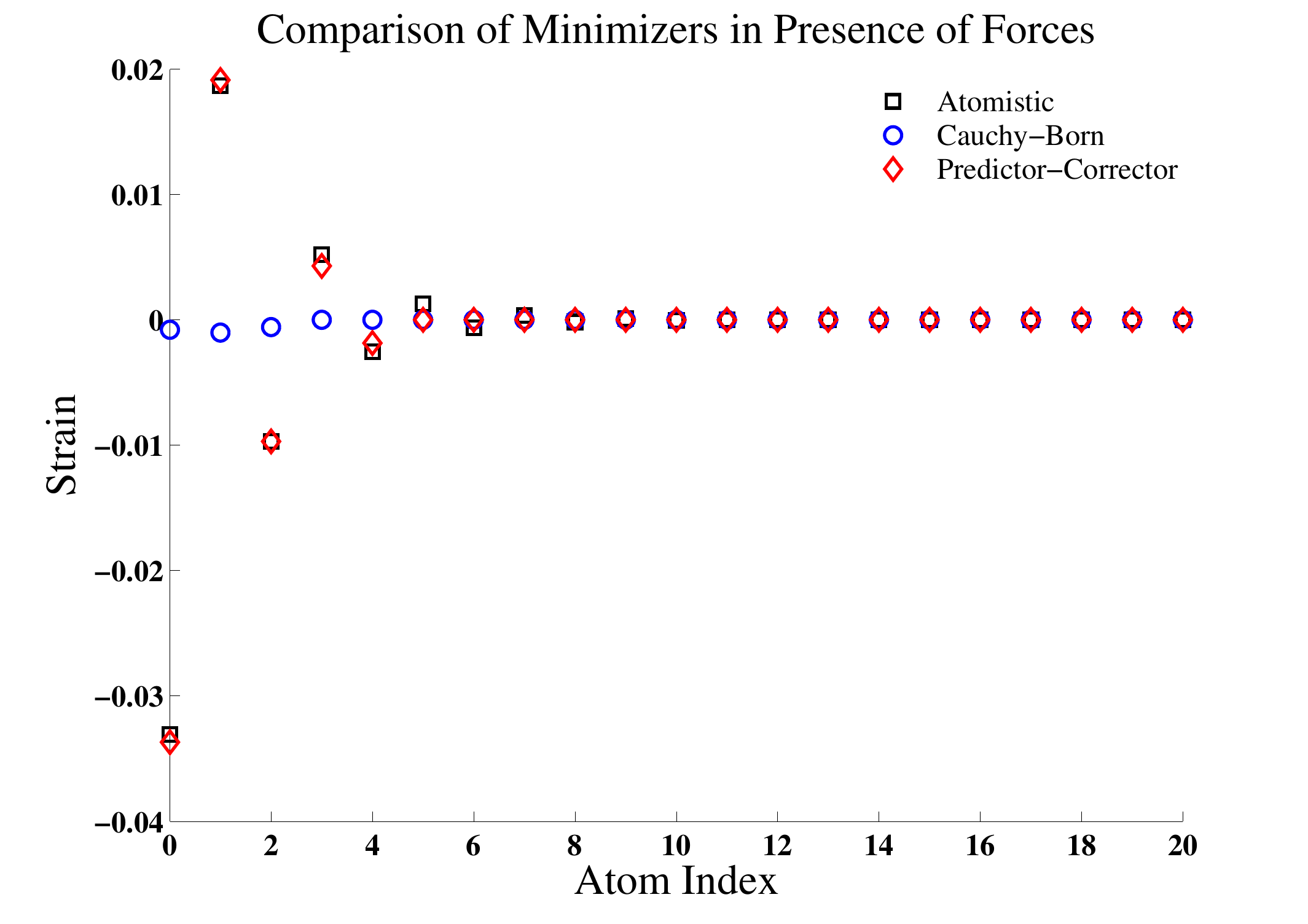}
     \caption{ }
     \label{fig:MinimizingStates}
   \end{subfigure}
   ~
   \begin{subfigure}[b]{0.48\textwidth}
    \includegraphics[scale=0.40]{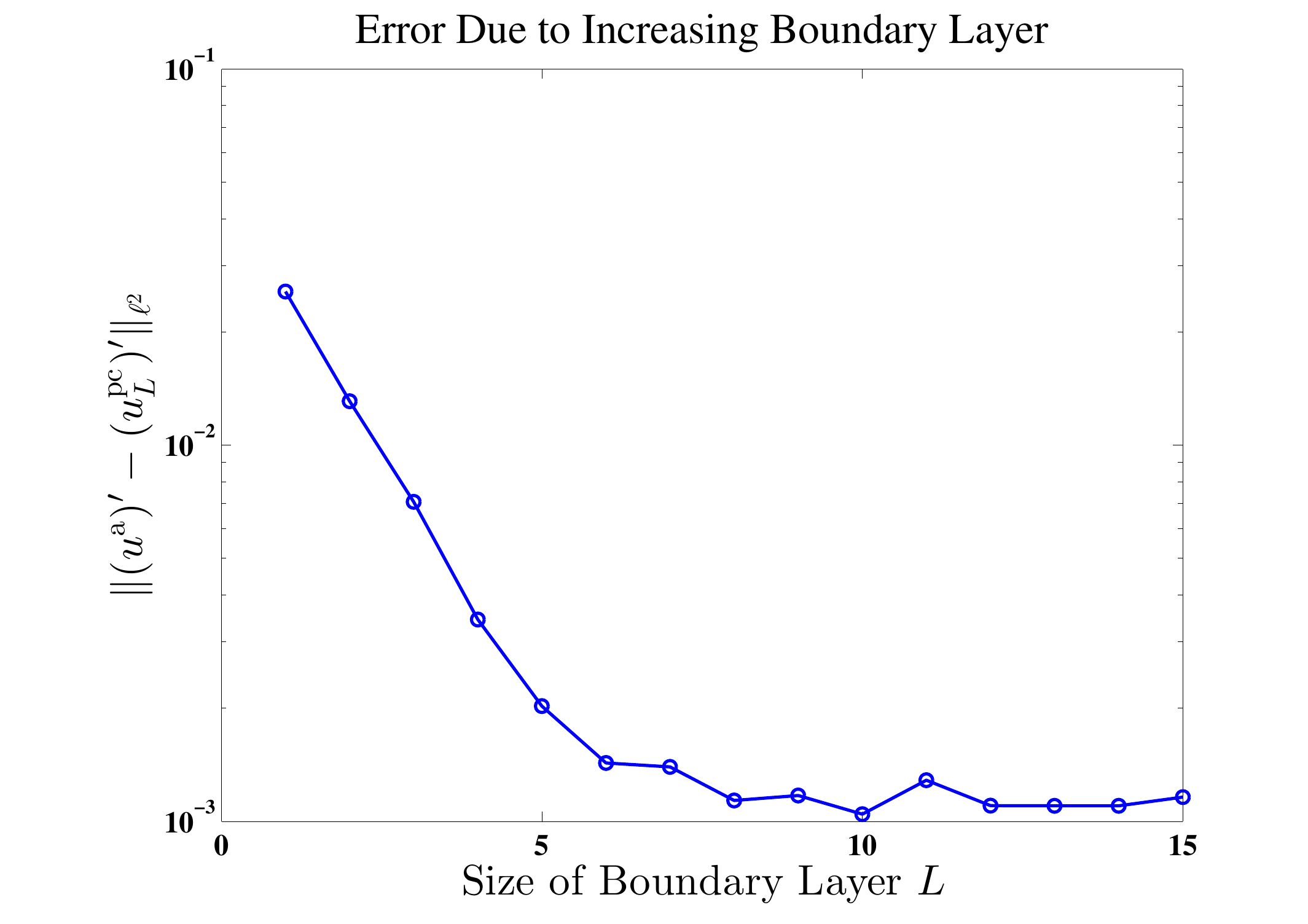}
    \caption{ }
    \label{fig:ResidualError}
   \end{subfigure}
   \caption{(a) Energy-minimizing configuration for atomistic and
     predictor-corrector systems in the presence of external forces.  Here,
     the size of the boundary layer for the predictor-corrector problem is $L
     = 5.$  (b) Demonstration of non-vanishing error in the presence of
     non-zero external forces and increasing boundary layer.}
\end{figure}

%
%
%
  We wish to highlight the fact that the error made by the predictor-corrector
  method's approximation will {\em not} vanish with increasing boundary layer in
  the case of {\em fixed} non-zero external forces.  To that end, we consider
  $f_{\ell} = \hat{f}(\ell)$, where $\hat{f}$ is given in \S~\ref{sec:test2}. In
  Figure \ref{fig:MinimizingStates}, we display the energy-minimizing
  configurations for the atomistic solution, the Cauchy--Born solution, and the
  predictor-corrector method with $L = 5$.  We observe that a small error in the
  boundary layer still persists.  Indeed, since the error depends on the
  magnitude and regularity of $f$ as demonstrated in
  \eqref{Eq:PCErrorInExtForce}, it cannot be driven to zero. This is numerically
  demonstrated by Figure \ref{fig:ResidualError}. We emphasize, however, that
  despite our extremely concentrated and non-smooth choice of $f$, the
  predictor-corrector method is able to reduce the error by a factor of $30$,
  which is an excellent improvement given the simplicity of the
  predictor-corrector scheme.


\subsection{Conclusions}
We produced a detailed analysis of an atomistic surface in a 1D model
problem. Based on this analysis, we showed that, while a regular Cauchy--Born
method fails to accurately capture the surface effects, it can be {\em
  corrected} (post-processed) using a computationally cheap surface
cell-problem. We gave a sharp error analysis of this new {\em
  predictor-corrector} scheme and, in particular, showed that its error relative
to the fully atomistic solution tends to zero in a natural scaling limit. Our
numerical results show promising quantitative behavior of the proposed scheme.

While the analysis is one-dimensional, we anticipate that it can shed some light
even on the three-dimensional case, when surface behavior dominates edge or
corner effects. An important new effect that will have to be taken into account
in two and three dimensions are surface stresses that act tangentially.



\section{Proofs: Atomistic Model}
\label{sec:proofs:atm}
%
%
%
\begin{proof}[Proof of Theorem \ref{Thm:ExistenceMinimizerAtomistic}]
  We employ the direct method of the calculus of variations. Property
  \eqref{Eq:MinSiteEnergy} on the site potential energies and the fact that
  $|V^{\surf}(0)| < \infty$ implies that $\E^{\rm a}(u) \geq 0$.  We may
  therefore consider an energy minimizing sequence $\{u_{n}\}_{n=1}^{\infty}
  \subset \U$ such that
  \begin{equation}\label{Eq:MinimizerProofAtomisticEnergyBound}
    \E^{\rm a}(u_{n}) \to \inf_{u \in \U}\E^{\rm a}(u)
    \qquad \text{and} \qquad
    \E^{\rm a}(u_{n}) < \inf_{u \in \U}\E^{\rm a}(u) + 1.
  \end{equation}

  We wish to show that $|u_{n}|_{H^{1}} = \|u_{n}'\|_{\ell^{2}(\Z_{\geq0})}$ is
  uniformly bounded. To that end, we will consider separately the contributions
  to the norm from the small and large strains.  Let $\eps > 0$.  For each
  $n \in \N$, we denote the set of indices that indicate the locations of
  $\eps$-defects in the chain to be $D^{n}_{\eps}$:
  \begin{displaymath}
    D^{n}_{\eps} :=
    \{\ell \in \Z_{\geq 0} \; | \; |u_{n,\ell}'| > \eps \; \text{or} \;
    |u_{n,\ell+1}'| > \eps \;\}.
  \end{displaymath}
  Property \eqref{Eq:PositiveSiteEnergy} implies that there exists a $\tau > 0$
  such that
  \begin{displaymath}
    \Vb(u_{n,\ell}', u_{n,\ell+1}') \geq \tau
    \qquad \text{for all} \; \ell \in D^{n}_\eps.
  \end{displaymath}
  Since $\Vb \geq 0$, we obtain
  \begin{displaymath}
    \E^{\rm a}(u_{n}) \geq
    \sum_{\ell \in D^{n}_{\eps}}\Vb(u'_{\ell}, u'_{\ell + 1})
    \geq
    \tau \#D^{n}_{\eps}.
  \end{displaymath}
  Hence, the cardinality of $D^{n}_{\eps}$ is uniformly bounded; that is,
  \begin{equation}\label{Eq:UniformDefectBound}
    \max_{n}\#D^{n}_{\eps} < \infty.
  \end{equation}

  We now analyze small strains.  By Taylor's Theorem and Property
  \eqref{Eq:PositiveDefiniteSiteEnergy}, we may show that for small enough
  $\eps$ there exists a constant $C_{\eps} > 0$ such that
  \begin{displaymath}
    \Vb(u_{n,\ell}', u_{n,\ell+1}')
    \geq
    C_{\eps}\big(|u_{n,\ell}'|^{2} + |u_{n,\ell+1}'|^{2}\big)
    \quad
    \text{for all} \quad \ell \notin D^{n}_{\eps}.
  \end{displaymath}
  With this inequality and the positivity of the site energies, we obtain
  \begin{align*}
 \E^{\rm a}(u_{n}) &=
   V^{\surf}(u_{n,0}')
     + \sum_{\ell \in \Z_{\geq 0}}
       \Vb(u_{n,\ell}', u_{n,\ell+1}')
 \\ &=
   V^{\surf}(u_{n,0}')
     + \sum_{\ell \in D^{n}_{\eps}}
       \Vb(u_{n,\ell}', u_{n,\ell+1}')
     + \sum_{\ell \notin D^{n}_{\eps}}
       \Vb(u_{n,\ell}', u_{n,\ell+1}')
  \\& \geq
      C_{\eps} \sum_{\ell \notin D^{n}_{\eps}}(|u_{n,\ell}'|^{2}
                                             + |u_{n,\ell+1}'|^{2}).
\end{align*}
From \eqref{Eq:MinimizerProofAtomisticEnergyBound}, we deduce the bound
\begin{equation}\label{Eq:MinimizerSmallStrainBound}
  \sum_{\ell \notin D^{n}_{\eps}}|u_{n,\ell}'|^{2}
  \leq
  C_{\eps}^{-1} \Bigl(\inf_{u \in \U}\E^{\rm a}(u) + 1\Bigr) < \infty.
\end{equation}

Note that the upper bound \eqref{Eq:MinimizerSmallStrainBound} excludes only a
finite number of strains on the chain (the $\eps$-defects).  If we can show that
there exists a positive constant $C > 0$ such that
\begin{equation}
  \label{eq:existence-groundstate-uniform-bound-u'}
  |u_{n,\ell}'| \leq C \qquad \text{for all} \quad n \in \N,
    \quad \ell \in \Z_{\geq 0}
\end{equation}
then we obtain a uniform bound on $\|u_n' \|_{\ell^2}$.

Suppose, for contradiction, that \eqref{eq:existence-groundstate-uniform-bound-u'}
fails. Then, there exists a subsequence $\{u_{n_{k}}\}_{k=1}^{\infty}$ and a
sequence of indices $\{\ell_{k}\}_{k=1}^{\infty}$ such that
\begin{equation}
  \label{Eq:MinimizerStrainToInfinity}
  \lim_{k \to \infty}|u_{n_{k},\ell_{k}}'| \to \infty.
\end{equation}
Since $\max_{n}\#D^{n}_{\eps}$ is finite, we may add the condition
on the indices $\ell_k$ that there exists a constant $S > 0$ such that
\begin{equation}
  \label{Eq:MinimizerAdjacentStrainBound}
  |u_{n_{k},\ell_{k} + 1}'| \leq S \quad \text{for all} \; k.
\end{equation}
We now split the chain into two components,
\begin{align*}
 \E^{\rm a}(u_{n_{k}})
  &=
    V^{\surf}(u_{n_{k},0}') +
    \sum_{\ell = 0}^{\ell_{k} - 1}\Vb(u_{n_{k},\ell}',
    u_{n_{k},\ell + 1}')
  \\ &+
       \Vb(u_{n_{k},\ell_{k}}', u_{n_{k},\ell_{k} + 1}') +
       \sum_{\ell = \ell_{k} + 1}^{\infty}\Vb(u_{n_{k},\ell}',
       u_{n_{k},\ell + 1}')  \\
  &=: E_1 + E_2.
\end{align*}
Since $\Vb \geq 0$, we obtain
\begin{displaymath}
  E_1 = V^{\surf}(u_{n_{k},0}') +
     \sum_{\ell = 0}^{\ell_{k} - 1}\Vb(u_{n_{k},\ell}',
                                       u_{n_{k},\ell + 1}')
 \geq
   \inf_{s' \in \R}V^{\surf}(s') > 0.
\end{displaymath}

To bound $E_2$, we observe that this group represents a semi-infinite chain that
has, up to a small error, decoupled from the first group. Therefore it is
bounded below by the infimum of the energy. To make this precise, let
$\eps_{1} := \frac{1}{2}\inf_{s \in \R}V^{\surf}(s) > 0$.  By Properties
\eqref{Eq:DefSurfSiteEnergy} and \eqref{Eq:SmoothnessSiteEnergy}, there exists
$R > 0$ such that if $|r| > R$, then $|V^{\bulk}(r, s) - V^{\surf}(s)| <
\eps_{1}$ for all $|s| \leq S$.  From \eqref{Eq:MinimizerStrainToInfinity},
there exists a $K > 0$ such that $|u_{n_{k},\ell_{k}}'| > R$ for all
$k \geq K$.  Therefore, we can conclude that, for $k \geq K$,
\begin{align*}
  E_2 &= \Vb(u_{n_{k},\ell_{k}}', u_{n_{k},\ell_{k} + 1}') +
        \sum_{\ell = \ell_{k} + 1}^{\infty}\Vb(u_{n_{k},\ell}',
        u_{n_{k},\ell + 1}')
 \\ &\geq
      V^{\surf}(u_{n_{k},\ell_{k} + 1}') - \eps_{1}
      + \sum_{\ell = \ell_{k} + 1}^{\infty}\Vb(u_{n_{k},\ell}',
      u_{n_{k},\ell + 1}')
  \\ &\geq
   \inf_{u \in \U} \E^{\rm a}(u) - \eps_{1}
\end{align*}
In summary, we have shown that, for $k \geq K$,
\begin{displaymath}
  \E^{\rm a}(u_{n_{k}}) = E_1 + E_2 \geq
  \smfrac{1}{2}\inf_{s \in \R}V^{\surf}(s)
  + \inf_{u \in \U} \E^{\rm a}(u).
\end{displaymath}
Since $\frac{1}{2}\inf_{s \in \R}V^{\surf}(s) > 0$, this contradicts the fact
that $\{u_{n_k} \}$ is a minimizing sequence.

Therefore, \eqref{eq:existence-groundstate-uniform-bound-u'} holds and together
with \eqref{Eq:UniformDefectBound} and \eqref{Eq:MinimizerSmallStrainBound},
and we deduce that $\| u_n' \|_{\ell^2}$ is uniformly bounded. Upon passing to a
subsequence (not relabeled), we may assume that there exists $u_\infty' \in
\ell^2(\Z_{\geq 0})$ such that $u_n' \rightharpoonup u_\infty'$ weakly in
$\ell^2$. In particular, this implies that
\begin{displaymath}
  u_{n,\ell}' \to u_{\infty,\ell}' \qquad \text{for all}
    \quad \ell \in \Z_{\geq 0}.
\end{displaymath}
We may now prove that $u_\infty$ is a minimizer by employing Fatou's Lemma to
estimate
\begin{align*}
 \inf_{u \in \U}\E^{\rm a}(u) &= \liminf_{n \to \infty}\E^{\rm a}(u_{n})
 \\ &\geq
      \liminf_{n \to \infty} V^{\surf}(u_{n,0}') +
      \sum_{\ell = 0}^{\infty}  \liminf_{n \to \infty}
      \Vb(u_{n,\ell}', u_{n,\ell + 1}')
= \Ea(u_\infty).
\end{align*}
Thus, $u_{\infty}$ is an energy minimizing configuration of the atomistic energy
$\E^{\rm a}$ in $\U$.
\end{proof}
%

%
%
\begin{proof}[Proof of Theorem \ref{Thm:ExpDecayAtomisticStrain}]
  Critical points of $\Ea$ satisfy the Euler--Lagrange equation
  \begin{displaymath}
    \langle \delta \E^{\rm a}(u'), v'\rangle
    =
    0
    \quad \text{for all} \; v \in \U.
  \end{displaymath}
  Expanding the first variation $\delta \E^{\rm a}(u')$ about $u' = 0$ yields
\begin{displaymath}
 \langle
   \delta \E^{\rm a}(0) + \delta^{2}\E^{\rm a}(0)u' + \zeta(u'), v'
 \rangle
 = 0,
\end{displaymath}
where $\zeta(u')$ is the remainder term from the expansion which can be readily
shown to satisfy the bounds
\begin{equation}
  \label{eq:exp-decay-bound-zeta}
  |\zeta_{0}| \lesssim |u_{0}'|^{2} + |u_{1}'|^{2} \qquad \text{and} \qquad
  |\zeta_{\ell}| \lesssim |u_{\ell-1}'|^{2} + |u_{\ell}'|^{2} + |u_{\ell+1}'|^{2}
  \quad \text{ for $\ell \geq 1$}.
\end{equation}
We may explicitly compute the remaining terms in the Taylor series expansion
keeping in mind that $\Vb(F,F)$ achieves its minimum at $(0,0)$:
\begin{align*}
 \langle \delta \E^{\rm a}(0), v' \rangle
   &=
 v_{0}'\Bigl\{\partial_{F}V^{\surf}(0)\Bigr\}
 \quad \text{and}
 \\
 \langle \delta^{2} \E^{\rm a}(0)u', v' \rangle
   &=
 v_{0}'u_{0}'\Bigl\{\partial^{2}_{F}V^{\surf}(0) +
   \partial^{2}_{11}\Vb(0, 0)\Bigr\}
 + v_{0}'u_{1}'\Bigl\{\partial^{2}_{12}\Vb(0,0)\Bigr\}
  \\
  + \sum_{\ell = 1}^{\infty}v_{\ell}'\biggl[
       &u_{\ell-1}'\Bigl\{\partial^{2}_{12}\Vb(0,0)\Bigr\}
      + u_{\ell}'\Bigl\{\partial^{2}_{22}\Vb(0,0) +
                        \partial^{2}_{11}\Vb(0,0)\Bigr\}
  \\ + &u_{\ell+1}'\Bigl\{\partial^{2}_{12}\Vb(0,0)\Bigr\}\biggr].
\end{align*}
Let
\begin{displaymath}
  b := \partial^{2}_{12}\Vb(0,0), \qquad
  a := \partial^{2}_{11}\Vb(0,0) + \partial^{2}_{22}\Vb(0,0),
  \quad \text{and} \quad
  a_{\rm s} := \partial^{2}_{F}V^{\surf}(0) +
   \partial^{2}_{11}\Vb(0,0),
\end{displaymath}
then we see that $u'$ solves the system
\begin{align}
  \label{Eq:SimpleDifferenceEquation-BC}
  a_{\rm s} u_{0}'
  +
  bu_{1}'
  &=
  -\zeta_{0} - \Bigl(\partial_{F}V^{\surf}(0)\Bigr)
  =:
  \tilde{\zeta}_{0}, \\
  %
\label{Eq:SimpleDifferenceEquation}
  bu_{\ell+1}'
  +
   a u_{\ell}'
  +
  b u_{\ell-1}'
  &=  -\zeta_{\ell} \qquad \text{for $\ell \geq 1$.}
\end{align}
%
%
%
%
We now suppose that $b \neq 0$.  The case when $b = 0$ is analogous.  Solving
for the homogeneous solution to this system of equations gives
\begin{displaymath}
  u_{\ell}' = c_{+}\lambda_{+}^{\ell} + c_{-}\lambda_{-}^{\ell},
  \qquad \text{where} \qquad \lambda_{\pm} = \frac{-a \pm \sqrt{a^{2} - 4b^{2}}}{2b}.
\end{displaymath}
Positive definiteness of $\nabla^2 \Vb(0,0)$ from Property
\eqref{Eq:PositiveDefiniteSiteEnergy} implies that $a - 2b > 0$ and
$a + 2b > 0$.  Hence, the discriminant is always positive, $a^{2} - 4b^{2} > 0$,
and as a result, $\lambda_{+} \neq \lambda_-$.  The symmetry of
\eqref{Eq:SimpleDifferenceEquation} implies that $\lambda_{+} = 1/\lambda_{-}$.
This relation combined with the fact that the discriminant is never 0 implies
that $\lambda_{\pm} \neq 1$.  Without loss of generality, then, we must have
that $0 < |\lambda_{+}| < 1$ and $|\lambda_{-}| > 1$.  As solutions to the
atomistic problem must belong to $\U$, we have that $u'_{\ell} \to 0$ as
$\ell \to \infty$.  This boundary condition at infinity implies that $c_{-} = 0$
in order to prevent exponential growth in the strain.  Thus, letting
$\lambda := \lambda_{+}$, we have that solutions of the homogeneous equation are
of the form
$ u_{\ell}' = u_0' \lambda^{\ell}$.
A discrete Green's function argument provides the solution for the inhomogeneous
equation:
\begin{displaymath}
  u_{\ell}' = C_{\text{BC}}\lambda^{\ell} +
  D^{-1} \sum_{k=0}^{\infty} \lambda^{|\ell-k|} \zeta_{k},
\end{displaymath}
where $D := \sqrt{a^{2} - 4b^{2}}$ and $C_{\rm BC}$ can be determined from
$u_0'$. Taking the absolute value of both sides, using the triangle inequality,
and applying \eqref{eq:exp-decay-bound-zeta}, we can estimate
\begin{align*}
|u_{\ell}'| &\leq C_{{\rm BC}}|\lambda|^{\ell}
  + D^{-1} \sum_{k=0}^{\infty} |\lambda|^{|\ell - k|} |\zeta_{k}|
\lesssim C_{\rm BC} |\lambda|^\ell
    + D^{-1} \sum_{k = 0}^{\infty}  |\lambda|^{|\ell - k|} |u_{k}'|^{2}.
\end{align*}
Let $\mu$ be a constant such that $|\lambda| < \mu < (1 + |\lambda|)/2$.  Then,
\begin{displaymath} 
\sum_{m = 0}^{\infty}\mu^{|\ell - m|}|u_{m}'| \lesssim
\sum_{m=0}^{\infty}C_{\text{BC}}\mu^{|\ell - m|}|\lambda|^{m}
  + D^{-1} \sum_{m,k \geq 0} \mu^{|\ell - m|}|\lambda|^{|m - k|} |u_{k}'|^{2}.
\end{displaymath}
Using the observation that
\begin{displaymath}
 \sum_{m = 0}^{\infty}\mu^{|\ell - m|}|\lambda|^{|m - k|} \lesssim \mu^{|\ell -
k|},
\end{displaymath}
we arrive at
%
\begin{equation}\label{Eq:InequalityUsedinExpDecayProof}
 \sum_{m = 0}^{\infty}\mu^{|\ell - m|}|u_{m}'| \lesssim C_{\text{BC}}\mu^{\ell}
+
 D^{-1} \sum_{k = 0}^{\infty} \mu^{|\ell - k|} |u_{k}'|^{2}.
\end{equation}
Ignoring the prefactor $D^{-1}$, the second term on the right-hand side of
\eqref{Eq:InequalityUsedinExpDecayProof} can be bounded by
\begin{displaymath}
 \sum_{k = 0}^{\infty} \mu^{|\ell - k|} |u_{k}'|^{2}
   =
 \sum_{k \leq k_{0}} \mu^{|\ell - k|} |u_{k}'|^{2}
   + \sum_{k \geq k_{0}} \mu^{|\ell - k|} |u_{k}'|^{2}
 \lesssim
 \mu^{\ell} + \sup_{k \geq k_{0}}|u_{k}'|\sum_{k \geq k_{0}} \mu^{|\ell -
k|} |u_{k}'|.
\end{displaymath}
For any $\eps > 0$, we can choose $k_0$ sufficiently large so that
$\sup_{k \geq k_{0}}|u_{k}'| \leq \eps$. Hence, we arrive at
\begin{displaymath}
 D^{-1}\sum_{k = 0}^{\infty}\mu^{|\ell - k|}|u_{k}'|^{2}
 \lesssim
 \mu^{\ell} + \eps\sum_{k \geq k_{0}}\mu^{|\ell - k|}|u_{k}'|.
\end{displaymath}
Substituting this bound into \eqref{Eq:InequalityUsedinExpDecayProof}
yields
\begin{displaymath}
  (1 - \eps)|u'_{\ell}|
  \leq
  (1 - \eps)\sum_{k=0}^{\infty}\mu^{|\ell - k|}|u_{k}'|
  \lesssim
  \mu^{\ell}. \qedhere
\end{displaymath}
%
\end{proof}
%
%
%

\section{Proofs: Cauchy--Born Model}
\label{sec:proofs-cb}
%
%
\begin{proof}[Proof of Proposition \ref{Thm:CBConsistencyErrorAtSurface}]
  Let $u = 0$.  Then, Properties
  \eqref{Eq:MinSiteEnergy}--\eqref{Eq:PositiveSiteEnergy} imply that
  $\delta \E^{\rm a}(u)$ is given by
  \begin{align*}
    \langle \delta \E^{\rm a}(u), v\rangle
    &=
      v_{0}'\biggl\{\partial_{F}V^{\surf}(u_{0}') +
      \partial_{1}\Vb(u_{0}',u_{1}')\biggr\}
      +
      \sum_{\ell = 1}^{\infty}v_{\ell}'\biggl\{
      \partial_{2}\Vb(u_{\ell - 1}',u_{\ell}')
      + \partial_{1}\Vb(u_{\ell}', u_{\ell + 1}')\biggr\} \\
    &= v_{0}' \partial_{F}V^{\surf}(0).
  \end{align*}
  This is maximized by taking $v_0' = {\rm sign}(\partial_{F}V^{\surf}(0))$ and
  $v_\ell' = 0$ for $\ell > 0$, thus proving the first statement.

  To deduce the lower bound on the error, we simply observe that
  \begin{displaymath}
    |\partial_{F}V^{\surf}(0)| =
    \| \delta\E^{\rm a}(0) \|_{\U^{*}} =
    \|\delta\E^{\rm a}(0) - \delta\E^{\rm a}(u^{\rm a}_{\rm gr})\|_{\U^{*}}
    \leq M \| (u^{\rm cb} - u^{\rm a}_{\rm gr})' \|_{\ell^2}. \qedhere
  \end{displaymath}
\end{proof}

Before we state the next result, recall that $f \in \U^*$ if and only there
exists $g \in \ell^2(\Z_{\geq 0})$ such that $\< f, u \> = \< g,
u'\>$. Summation by parts, with the convention  $g_{-1} = 0$,
may be used to show that
\begin{equation}
  \label{eq:f=devg-discrete}
  f_{\ell} = g_{\ell} - g_{\ell - 1}.
\end{equation}
Conversely, if we are given $f$, we may recover $g$ via
\begin{equation}
  \label{eq:g=sumf-discrete}
  g_\ell := \sum_{k = \ell + 1}^\infty f_k
  := \lim_{K \to \infty} \sum_{k = \ell + 1}^K f_k.
\end{equation}

In the Cauchy--Born model, we identify $f$ with its piecewise affine
interpolant. The continuous analogue of \eqref{eq:g=sumf-discrete} is
\begin{equation}
  \label{eq:gtil=intf-cts}
  \tilde{g}(x) := \int_x^\infty f(s) \, ds :=
  \lim_{K \to \infty} \int_x^K f(s) \, ds.
\end{equation}

\begin{lemma}
  Let $f \in \U^*$. Then, $g$ and $\tilde{g}$ are well-defined and satisfy the
  estimates
  \begin{align}
    \label{eq:estimate-g-gtil}
    \big| g_\ell - \tilde{g}(\ell+1/2) \big| &\leq \smfrac18 |f_\ell'|
                                               \qquad \text{for all}
                                               \quad \ell \in \Z_{\geq 0},
    \qquad \text{and} \\
    \label{eq:bound-gtil-lt-g}
    |\tilde{g}(x)| &\leq C \sum_{m = \ell-1}^{\ell+1} |g_m| \qquad
                     \text{for all} \quad x \in [\ell, \ell+1].
  \end{align}
\end{lemma}
\begin{proof}
  Since $f \in \U^*$ there exists $g \in \ell^2$ such that
  \eqref{eq:f=devg-discrete} holds. Let $M > \ell+1$.  Then,
  \begin{displaymath}
    \sum_{m = \ell+1}^M f_m = \sum_{m = \ell+1}^M (g_{\ell-1} - g_\ell)
    = g_\ell - g_M \to g_\ell \qquad \text{as $M \to \infty$.}.
  \end{displaymath}
  Hence, \eqref{eq:g=sumf-discrete} is well-defined.

  Similarly, let $M \in \Z$ and $M > \ell+1$.  Then,
  \begin{align*}
    \int_{\ell+1/2}^M f(s) \,ds  - \sum_{m = \ell+1}^M f_m
    &= \int_{\ell+1/2}^{\ell+1} f(s) \,ds
      + \int_{\ell+1}^M f(s) \, ds - \sum_{m = \ell+1}^M f_m \\
    &= \smfrac12 \big(\smfrac14 (f_\ell+f_{\ell+1}) + \smfrac12 f_{\ell+1} \big)
      + \sum_{m = \ell+1}^{M-1} \Big[ \smfrac12 (f_m + f_{m+1})
      - f_m \Big] - f_M \\
    &= \smfrac18 (f_\ell - f_{\ell+1}) - \smfrac12 f_{M}
      \to \smfrac18 (f_\ell - f_{\ell+1}) \quad \text{as $M \to \infty$.}
  \end{align*}
  Thus, we can conclude that $\tilde{g}(x)$ is well-defined for all $x$, as well
  as the stated estimate \eqref{eq:estimate-g-gtil}.

  To prove the second estimate, we simply note that, for $x \in [\ell, \ell+1]$,
  \begin{displaymath}
    |\tilde{g}(x) - \tilde{g}(\ell+1/2)| = \bigg| \int_{\ell+1/2}^x f(s) \,ds \bigg|
    \lesssim |f_\ell| + |f_{\ell+1}|,
  \end{displaymath}
  and then apply \eqref{eq:f=devg-discrete}.
\end{proof}

Next, we prove Theorem \ref{Thm:CBExistenceRegularity}.  In addition to the
bounds stated therein, we will also
establish that
\begin{equation}
  \label{eq:CBExistenceRegularity-Bound-Du}
  |\nabla u^{\rm cb}(x)| \lesssim |\tilde{g}(x)|,
\end{equation}

%
%
\begin{proof}[Proof of Theorem \ref{Thm:CBExistenceRegularity} and of
  \eqref{eq:CBExistenceRegularity-Bound-Du}]
  The Euler--Lagrange equation of the Cauchy--Born problem
  \eqref{SemiInfiniteCauchyBornProblem} is
\begin{equation}\label{Eq:ELinCBExistenceProof}
 -\nabla \bigl[\partial_{F}W(\nabla u)\bigr] = f.
\end{equation}
%
%
Formally integrating over $(x, \infty)$ and using the fact that
$\partial_F W(0) = 0$, we obtain
\begin{displaymath}
  \partial_{F} W(\nabla u) = \int_{x}^{\infty}f(s)ds = \tilde{g}(x).
\end{displaymath}
We will now produce a function $u \in H^3_{\rm loc}$ satisfying this equation.

%

From Properties \eqref{Eq:PositiveDefiniteSiteEnergy} and
\eqref{Eq:SmoothnessSiteEnergy}, there exists $G > 0$ such that
$\partial^{2}_{F}W(F) > \frac{1}{2}\partial_{F}^{2}W(0) > 0$ on $[-G,G]$.  Thus,
$\partial_{F}W(F)$ is strictly monotone on this interval. The bound
\eqref{eq:bound-gtil-lt-g} implies that, for $\varepsilon_{\rm cb}$ sufficiently
small, $\tilde{g}([0,\infty)) \subset \partial_{F}W([-G,G])$.  Hence, we can define
\begin{displaymath}
\nabla u(x) = (\partial_{F}W)^{-1}(\tilde{g}(x)),
\end{displaymath}
where $(\partial_{F}W)^{-1}$ is defined to have range in $[-G, G]$.
%
Due to the restriction on the range of the inverse function,
$\nabla u(x) \in [-G,G]$, hence it follows that the solution must be stable with
a stability constant $c_{\rm cb} \geq \frac{1}{2}\partial^{2}_{F}W(0)$.

It is further easy to show with this that $(\partial_{F}W)^{-1}$ is twice
continuously differentiable, and in particular, Lipschitz. Noting that
$(\partial_{F}W)^{-1}(0) = 0$, Lipschitz continuity then yields the bound
\begin{displaymath}
  |\nabla u^{\rm cb}(x)|
  =
  |(\partial_{F}W)^{-1}(\tilde{g}(x)) - (\partial_{F}W)^{-1}(0)|
  \lesssim
  |\tilde{g}(x)|,
\end{displaymath}
where we have ignored the Lipschitz constant for the inverse function.

The remaining estimates are consequences of the fact that
$\tilde{g} \in C^{1,1}(\R_{\geq 0})$ with $\nabla \tilde{g} = -f$ and
$\nabla^2 \tilde{g} = -\nabla f$ along with an elementary computation.
\end{proof}

Next, we recall an auxiliary result that allows us to reduce the continuous
Cauchy--Born model to a discrete model. To that end, we recall that we can
identify discrete test functions $v \in \U$ with continuous test functions
$v \in \U^{\rm cb}$ via piecewise affine interpolation. Through the same
identification, we can also admit $u \in \U$ as arguments for $\del\E^{\rm cb}$.

\begin{proposition}\label{Proof:CBDiscretization}
  Under the conditions of Theorem \ref{Thm:CBExistenceRegularity}, we have
  %
  \begin{displaymath}
    \big|\big\langle \delta \E^{\rm cb}(u^{\rm cb})
                   - \delta \E^{\rm cb}(\Pi_{\rm a}u^{\rm cb}), v
    \big\rangle\big|
    \lesssim \|\nabla^{2}u^{\rm cb}\|_{L^{4}}^2  \|\nabla v\|_{L^{2}}
    \qquad \text{for all} \quad v \in \U.
  \end{displaymath}
\end{proposition}
\begin{proof}
  This result is a simplified variant of \cite[Lemma 5.2]{acta.atc}.
\end{proof}

\section{Proofs: Corrector Problem}
\label{sec:proofs-corrector}
%
%
\begin{proof}[Proof of Theorem \ref{Thm:ExistenceInfiniteCorrector}]
  {\it Consistency: } It is reasonable to expect that, for $|F_0|$ small, the
  corrector solution $q_\infty$ is close to $u^{\rm a}_{\rm gr}$, so we may find
  $q_\infty$ be applying the inverse function theorem in a neighborhood of
  $u^{\rm a}_{\rm gr}$. From Proposition \ref{Prop:FirstVariations}, the first
  variation of the corrector energy evaluated at
  $\tilde{q} = u^{\rm a}_{\rm gr}$ is
  \begin{align*}
    \langle \delta \E^{\Gamma}(\tilde{q}; F_0), v \rangle
    = v_{0}'\biggl\{\partial_{F}V^{\surf}(F_{0} + \tilde{q}_{0}') +
                    \partial_{1}\Vb(F_{0}+\tilde{q}_{0}', F_{0}+\tilde{q}_{1}')-
                    \partial_{F}W(F_{0})\biggr\}&
    \\ +
         \sum_{\ell = 1}^{\infty}v_{\ell}'\biggl\{
         \partial_{2}\Vb(F_{0} + \tilde{q}_{\ell-1}', F_{0} + \tilde{q}_{\ell}')
         +\partial_{1}\Vb(F_{0} + \tilde{q}_{\ell}', F_{0} +
\tilde{q}_{\ell+1}')
         -\partial_{F}W(F_{0})  \biggr\} &\,,
  \end{align*}
  where $v \in \Q_{\infty}$. Using the fact that $\del\Ea(\tilde{q}) = 0$ and
  $\partial_F W(0) = 0$, we obtain
  \begin{align*}
    \langle \delta \E^{\Gamma}(\tilde{q}; F_0), v \rangle
    &= \langle \delta \E^{\Gamma}(\tilde{q}; F_0) - \del\Ea(\tilde{q}), v \rangle  \\
    &=   v_{0}'\biggl\{
      \partial_{F}V^{\surf}(F_{0} + \tilde{q}_{0}')
      + \partial_{1}\Vb(F_{0}+\tilde{q}_{0}', F_{0}+\tilde{q}_{1}')
      - \partial_{F}W(F_{0})  \\
    & \qquad \qquad
      - \partial_{F}V^{\surf}(\tilde{q}_{0}')
      - \partial_{1} \Vb( \tilde{q}_{0}', \tilde{q}_{1}')
      + \partial_{F}W(0) \biggr\}  \\
    & \quad +
      \sum_{\ell = 1}^{\infty}v_{\ell}'\biggl\{
      \partial_{2}\Vb(F_{0} + \tilde{q}_{\ell-1}', F_{0} + \tilde{q}_{\ell}')
      + \partial_{1}\Vb(F_{0} + \tilde{q}_{\ell}', F_{0} + \tilde{q}_{\ell+1}')
      - \partial_{F}W(F_{0})   \\
    & \qquad \qquad \qquad
      - \partial_{2}\Vb( \tilde{q}_{\ell-1}', \tilde{q}_{\ell}')
      - \partial_{1}\Vb( \tilde{q}_{\ell}',  \tilde{q}_{\ell+1}')
      + \partial_{F}W(0)
      \biggr\} \\
    &=: \sum_{\ell = 0}^\infty v_\ell' A_\ell.
  \end{align*}
  Using Lipschitz continuity of all potential derivatives we easily obtain that
  \begin{align*}
    |A_0| &\lesssim |F_0|.
  \end{align*}
  To estimate $A_\ell$ for $\ell \geq 1$, we proceed more carefully. Expanding with
  respect to $\tilde{q}_j$ for $j = \ell-1, \ell, \ell+1$, employing the identity
  $\partial_F W(F) = \partial_1 \Vb(F, F) + \partial_2 \Vb(F, F)$, and the
  global Lipschitz continuity of $\partial^2 \Vb$, we obtain
  \begin{align*}
    |A_\ell|
    &= \bigg| \Big\{\partial_2 \Vb(F_0, F_0) + \partial_1 \Vb(F_0, F_0)
      - \partial_F W(F_0)\Big\}
      - \Big\{ \partial_2 \Vb(0, 0) + \partial_1 \Vb(0, 0) - \partial_F W(0) \Big\} \\
    & \qquad \qquad
      + \int_0^1 \Big\{\partial \partial_2
      \Vb(F_0 + t \tilde{q}_{\ell-1}', F_0 + t\tilde{q}_{\ell}')
      - \partial \partial_2 \Vb(t \tilde{q}_{\ell-1}',  t\tilde{q}_{\ell}') \Big\}
      \cdot
      \begin{pmatrix}
        \tilde{q}_{\ell-1}' \\
        \tilde{q}_\ell'
      \end{pmatrix} \, dt \\
    & \qquad \qquad
      + \int_0^1 \Big\{\partial \partial_1
      \Vb(F_0 + t \tilde{q}_{\ell}', F_0 + t\tilde{q}_{\ell+1}')
      - \partial \partial_1 \Vb(t \tilde{q}_{\ell}',  t\tilde{q}_{\ell+1}') \Big\}
      \cdot
      \begin{pmatrix}
        \tilde{q}_{\ell}' \\
        \tilde{q}_{\ell+1}'
      \end{pmatrix}
    \, dt
    \bigg| \\
    &\lesssim |F_0| \big( |\tilde{q}_{\ell-1}'|
      + |\tilde{q}_\ell'| + |\tilde{q}_{\ell+1}'| \big).
  \end{align*}
  An application of the Cauchy--Schwartz inequality yields
  \begin{equation}
    \label{eq:proof-existence-corrector:consistency}
     \langle \delta \E^{\Gamma}(\tilde{q}; F_0), v \rangle
     \leq \|A \|_{\ell^2} \|v' \|_{\ell^2}
     \lesssim |F_0| (1 + \|\tilde{q}'\|_{\ell^2}) \| v' \|_{\ell^2}.
  \end{equation}

  {\it Stability: } Recall that $\tilde{q} = u^{\rm a}_{\rm gr}$ is strongly
  stable in the atomistic model \eqref{Eq:AtomisticGroundStateStability} with
  stability constant $c_{\rm a} > 0$. To prove stability of
  $\ddel\E^{\Gamma}(\tilde{q}; F_0)$, we simply bound the error in the Hessians:
  { \small
  \begin{align*}
    &\big\< \big(\ddel\Ea(\tilde{q})
      - \ddel\E^\Gamma(\tilde{q};F_0)\big) v, v \big\>
    = |v_{0}'|^{2} \Bigl\{\partial^{2}_{F}V^{\surf}(\tilde{q}_{0}')
      - \partial^{2}_{F}V^{\surf}(F_0 + \tilde{q}_{0}')  \Bigr\} \\
    & +
      \sum_{\ell = 0}^{\infty}\biggl[
      |v_{\ell}' + v_{\ell + 1}'|^{2}
      \Bigl\{\partial^{2}_{12}\Vb(\tilde{q}_{\ell}', \tilde{q}_{\ell + 1}')
      - \partial^{2}_{12}\Vb(F_0 + \tilde{q}_{\ell}', F_0+ \tilde{q}_{\ell +
1}')
      \Bigr\} \\
    &\quad \qquad
      + |v_{\ell}'|^{2}
      \Bigl\{ \partial^{2}_{11}\Vb(\tilde{q}_{\ell}',\tilde{q}_{\ell+1}')
         - \partial^{2}_{12}\Vb(\tilde{q}_{\ell}',\tilde{q}_{\ell+1}')
      -\partial^{2}_{11}\Vb(F_0 + \tilde{q}_{\ell}', F_0 + \tilde{q}_{\ell+1}')
         + \partial^{2}_{12}\Vb(F_0 + \tilde{q}_{\ell}',F_0 +
\tilde{q}_{\ell+1}')
      \Bigr\} \\
    & \quad\qquad
      + |v_{\ell+1}'|^{2}
      \Bigl\{\partial^{2}_{22}\Vb(\tilde{q}_{\ell}',\tilde{q}_{\ell+1}')
      - \partial^{2}_{12}\Vb(\tilde{q}_{\ell}',\tilde{q}_{\ell+1}')
      - \partial^{2}_{22}\Vb(F_0 + \tilde{q}_{\ell}',F_0 + \tilde{q}_{\ell+1}')
      + \partial^{2}_{12}\Vb(F_0 + \tilde{q}_{\ell}',F_0 + \tilde{q}_{\ell+1}')
      \Bigr\} \biggr].
  \end{align*}
}Employing Lipschitz continuity of $\partial_F^2 \Vs$ and $\partial^2 \Vb$, we
can immediately deduce that
\begin{align*}
  \big| \big\< \big(\ddel\Ea(\tilde{q})
      - \ddel\E^\Gamma(\tilde{q};F_0)\big) v, v \big\>  \big|
  \leq  C |F_0| \| v' \|_{\ell^2}^2,
\end{align*}
where the constant $C$ is the upper bound on the Lipschitz constants
involved in the estimate multiplied by a simple factor. Thus, we obtain, for
$|F_0| \leq \frac14 c_{\rm a}/C$,
\begin{equation}
  \label{eq:proof-existence-corrector:stab}
  \langle \delta^{2} \E^{\Gamma}(\tilde{q}; F_0)v, v \rangle
  \geq
  \frac{3 c_{\rm a}}{4}\|v'\|_{\ell^{2}}^{2} \qquad
    \text{for all} \quad v \in \Q_\infty.
\end{equation}

The consistency estimate \eqref{eq:proof-existence-corrector:consistency}, the
stability estimate \eqref{eq:proof-existence-corrector:stab}, Lipschitz
continuity of $q \mapsto \delta^{2} \E^{\Gamma}(q; F_0)$, and an
application of Theorem \ref{Thm:InverseFunctionTheorem} implies the existence of
an equilibrium $q_\infty$ with
$\| \tilde{q}' - q_\infty' \|_{\ell^2} \lesssim |F_0|$. Choosing $F_0$
sufficiently small and once again using the Lipschitz continuity of
$\delta^{2} \E^{\Gamma}$ implies \eqref{eq:corrector-stab}.

The exponential decay in \eqref{eq:corrector-exp-decay} follows from a
straightforward modification of Theorem \ref{Thm:ExpDecayAtomisticStrain}.
\end{proof}

To analyze the projection of the corrector problem from $\Q_\infty$ to $\Q_L$,
we introduce a projection operator $\Pi_L : \Q_{\infty} \to \Q_L$ via
\begin{displaymath}
 \bigl(\Pi_{L}q_{\ell}\bigr)' :=
   \left\{
     \begin{array}{lr}
       q_{\ell}', & \ell = 0, \ldots, L, \\
       0, & \ell > L.
     \end{array}
   \right.
\end{displaymath}

\begin{lemma}\label{lemma:ExpDecayTruncation}
  Let $q \in \Q_{\infty}$ satisfy $|q_\ell'| \lesssim \mu^\ell$ for some
  $\mu \in [0, 1)$.  Then,
\begin{displaymath}
 \|q' - \Pi_{L}q'\|_{\ell^{2}} \lesssim \mu^{L}.
\end{displaymath}
\end{lemma}
\begin{proof}
First, note that
\begin{displaymath}
 q_{\ell}' - \bigl(\Pi_{L}q\bigr)_{\ell}' =
   \left\{
     \begin{array}{lr}
       0, & 0 \leq \ell \leq L - 1, \\
       q_{\ell}', & \ell > L.
     \end{array}
   \right.
\end{displaymath}
Thus,
\begin{displaymath}
 \|q' - \Pi_{L}q'\|^{2}_{\ell^{2}} =
  \sum_{\ell = L}^{\infty}|q_{\ell}'|^{2}
  \lesssim
  \sum_{\ell = L}^{\infty}\mu^{2\ell} \\
  =
  \frac{\mu^{2L}}{1 - \mu^{2}} \qedhere
\end{displaymath}
\end{proof}
%

%
\begin{proof}[Proof of Proposition \ref{Prop:ExistenceFiniteCorrector}]
  The result is proven by an application of the quantitative inverse function
  theorem, Theorem \ref{Thm:InverseFunctionTheorem}, using the projected
  solution $\Pi_{L}q_{\infty}$ as an approximate solution. Observe that
  \begin{displaymath}
    \langle \delta \E^{\Gamma}(q_{\infty}; F_0), v \rangle = 0
    \qquad \text{for all} \quad v \in \Q_L \subset \Q_\infty,
  \end{displaymath}
  where $q_{\infty}$ solves the infinite corrector problem. By the Lipschitz
  continuity of $\delta \E^\Gamma$,
  \begin{displaymath}
    |\langle \delta \E^{\Gamma}(\Pi_{L}q_{\infty}; F_0), v \rangle|
    =
    |\langle \delta \E^{\Gamma}(\Pi_{L}q_{\infty}; F_0)
    - \delta \E^{\Gamma}(q_{\infty}; F_0), v \rangle|
    \lesssim
    \|q_{\infty}' - \Pi_{L}q_{\infty}'\|_{\ell^{2}}\|v'\|_{\ell^{2}}.
  \end{displaymath}
  According to Theorem \ref{Thm:ExistenceInfiniteCorrector}, there exists
  $0 \leq \mu_{q} < 1$ such that
  $|(q_{\infty})_{\ell}'| \lesssim \mu_{q}^{\ell}$.  Hence, Lemma
  \ref{lemma:ExpDecayTruncation} may be applied to show that
  \begin{displaymath}\label{Eq:CorrectorOperatorNormInequality}
    \|\delta  \E^{\Gamma}(\Pi_{L}q_{\infty}; F_0)\|_{\Q_{L}^{*}}
    = \sup_{v \in \Q_{L},\|v'\|_{\ell^{2}}=1}|\langle \delta
    \E^{\Gamma}(\Pi_{L}q_{\infty}; F_0), v \rangle| \lesssim \mu_{q}^{L}.
  \end{displaymath}

  Stability of $\Pi_{L}q_{\infty}$ follows from the stability of $q_{\infty}$
  (Theorem \ref{Thm:ExistenceInfiniteCorrector}), the Lipschitz continuity of
  $\ddel\E^\Gamma$, and Lemma \ref{lemma:ExpDecayTruncation},
\begin{align*}
 \langle \delta \E^{\Gamma}(\Pi_{L}q_{\infty}; F_0)v, v \rangle
 &=
 \langle \delta \E^{\Gamma}(q_{\infty}; F_0)v, v \rangle
 - \biggl[
     \langle \delta \E^{\Gamma}(\Pi_{L}q_{\infty}; F_0)v, v \rangle
   - \langle \delta \E^{\Gamma}(q_{\infty}; F_0)v, v \rangle
   \biggr]
 \\ &\geq
   \langle \delta \E^{\Gamma}(q_{\infty}; F_0)v, v \rangle
   -
   M \|q_{\infty}' - \Pi_{L}q_{\infty}'\|_{\ell^{2}}\|v'\|^{2}_{\ell^{2}}
 \\ &\geq
      \left(\frac{c_{\rm a}}{2} - CM \mu_{q}^{L}\right)\|v'\|_{\ell^{2}},
\end{align*}
where $M$ is the Lipschitz constant for $\E^\Gamma$ and $C$ is the
unlisted constant in \eqref{Eq:CorrectorOperatorNormInequality}.

For $L$ sufficiently large, all assumptions of Theorem
\ref{Thm:InverseFunctionTheorem} are met, and its application
completes the proof.
\end{proof}

\section{Proofs: Predictor-Corrector Method}
\label{sec:proofs-pc}
The center-piece of our analysis of the predictor-corrector method is the
following consistency error estimate. In its statement, we employ the notation
\begin{displaymath}
  u_{\ell}'' := u_{\ell + 1}' - u_{\ell}' \qquad \text{and} \qquad
  u_{\ell}''' := u_{\ell - 1}' - 2u_{\ell}' + u_{\ell}'.
\end{displaymath}

\begin{theorem}\label{Thm:SingleElementConsistencyError}
  Let $w, q \in \U$, $u := w + q$, and $F_{0} \in \R$.  Then, for all
  $v \in \U$,
 \begin{align*}
   |\langle \delta \E^{\rm a}(u) - \delta \E^{\rm cb}(w) - \delta
\E^{\Gamma}(q; F_0), v
    \rangle|
   \lesssim
   \Biggl[ \|w''\|_{\ell^{4}(\Z_{\geq 0})}^{2} +
         \|w'''\|_{\ell^{2}(\Z_{>0})} +
         \|q'' \cdot w''\|_{\ell^{2}(\Z_{\geq 0})}&  \\
        +
         |w_1' - F_0| + \Biggl(\sum_{\ell = 1}^{\infty}
           (|q_{\ell - 1}'| + |q_{\ell}'| + |q_{\ell + 1}'|)^{2}
            |w_{\ell}' - F_0|^{2}\Biggr)^{1/2}  &\,\Biggr]
   \cdot \|v'\|_{\ell^{2}}.
 \end{align*}
\end{theorem}
\begin{proof}
The difference in the first variations is given by
\begin{align*}
 \langle \delta &\E^{\rm a}(u) - \delta \E^{\rm cb}(w)
   - \delta \E^{\Gamma}(q; F_{0}), v
 \rangle
 \\ &=
 v_{0}'\biggl[
   \Bigl\{\partial_{F}V^{\surf}(u_{0}') +
          \partial_{1}\Vb(u_{0}', u_{1}')\Bigr\}
 - \Bigl\{\partial_{F}W(F_{0})\Bigr\}
 \\ & \qquad \qquad - \Bigl\{\partial_{F}V^{\surf}(F_{0} + q_{0}') +
              \partial_{1}\Vb(F_{0} + q_{0}', F_{0} + q_{1}') -
              \partial_{F}W(F_{0})\Bigr\}\biggr]
 \\ &\qquad +
 \sum_{\ell = 1}^{\infty}
   v_{\ell}'\biggl[
    \Bigl\{\partial_{2}\Vb(u_{\ell - 1}', u_{\ell}') +
           \partial_{1}\Vb(u_{\ell}', u_{\ell + 1}')\Bigr\}
      - \Bigl\{\partial_{F}W(w_{\ell}') \Bigr\}
  \\ &\qquad \qquad \qquad- \Bigl\{
    \partial_{2}\Vb(F_{0} + q_{\ell - 1}', F_{0} + q_{\ell}') +
    \partial_{1}\Vb(F_{0} + q_{\ell}', F_{0} + q_{\ell + 1}') -
    \partial_{F}W(F_{0})\Bigr\}\biggr] \\
  &=: S + B,
\end{align*}
where $S, B$ denote, respectively, the surface and bulk contributions to the
consistency error.

{\it Surface term: } Using the Lipschitz continuity of the bulk site energy, we
can show that
\begin{align}
  |S| &=  \biggl|v_{0}' \Bigl[
   \Bigl\{\partial_{F}V^{\surf}(u_{0}') +
          \partial_{1}\Vb(u_{0}', u_{1}')\Bigr\}
 - \Bigl\{\partial_{F}W(F_{0})\Bigr\}
  \notag
  \\ & \qquad \qquad -
   \Bigl\{\partial_{F}V^{\surf}(F_{0} + q_{0}') +
          \partial_{1}\Vb(F_{0} + q_{0}', F_{0} + q_{1}') -
          \partial_{F}W(F_{0})\Bigr\}\Bigr]\biggl|
  \notag
 \\ &= |v_0'| \,
      \Bigl|
   \partial_{1}\Vb(F_{0} + q_{0}', w_{1}' + q_{1}')
 - \partial_{1}\Vb(F_{0} + q_{0}', F_{0} + q_{1}')\Bigr|
  \notag
 \\ &\lesssim
 |v_{0}'| \cdot |w_{1}' - F_{0}| .
       \label{eq:pc:mainlemma:estiamte-S}
\end{align}

{\it Bulk sum: } We write $B = \sum_{\ell = 1}^\infty v_\ell' B_\ell$, where
\begin{align*}
  B_\ell &= \Bigl\{\partial_{2}\Vb(u_{\ell - 1}', u_{\ell}') +
           \partial_{1}\Vb(u_{\ell}', u_{\ell + 1}')\Bigr\}
      - \Bigl\{\partial_{F}W(w_{\ell}') \Bigr\}
  \\ &\qquad \qquad \qquad- \Bigl\{
    \partial_{2}\Vb(F_{0} + q_{\ell - 1}', F_{0} + q_{\ell}') +
    \partial_{1}\Vb(F_{0} + q_{\ell}', F_{0} + q_{\ell + 1}') -
    \partial_{F}W(F_{0})\Bigr\}.
\end{align*}
Using the identity
$\partial_F W(F) = \partial_1 \Vb(F, F) + \partial_2 \Vb(F, F)$ and adding 0,
we can split $B_\ell$ into {\small
\begin{align*}
B_\ell  &=\biggl\{\Bigl[\partial_{2}\Vb(u_{\ell-1}', u_{\ell}')
       + \partial_{1}\Vb(u_{\ell}', u_{\ell + 1}')\Bigr]
  -
 \Bigl[\partial_{2}\Vb(w_{\ell}' + q_{\ell -1}', w_{\ell}' + q_{\ell}')
      + \partial_{1}\Vb(w_{\ell}' + q_{\ell}', w_{\ell}' + q_{\ell + 1}')
 \Bigr]\biggr\}
 \\ & \qquad +
      \biggl\{\Bigl[\partial_{2}\Vb(w_{\ell}' + q_{\ell -1}', w_{\ell}' +
q_{\ell}')
      + \partial_{1}\Vb(w_{\ell}' + q_{\ell}', w_{\ell}' + q_{\ell + 1}')
 \Bigr]
 -
 \Bigl[\partial_{2}\Vb(w_{\ell}', w_{\ell}')
      + \partial_{1}\Vb(w_{\ell}', w_{\ell}')
 \Bigr] \biggr\}
 \\ & \qquad -
 \biggl\{\partial_{2}\Vb(F_{0} + q_{\ell - 1}', F_{0} + q_{\ell}')
      + \partial_{1}\Vb(F_{0} + q_{\ell}', F_{0} + q_{\ell + 1}')
      -  \partial_{2}\Vb(F_{0}, F_{0})
      - \partial_{1}\Vb(F_{0}, F_{0})
 \biggr\} \\
  &=: B_\ell^{(1)} + B_\ell^{(2)} - B_\ell^{(3)}.
\end{align*}}

{\it Estimate for $B_\ell^{(1)}$: } Expanding $\partial \Vb$, we obtain, using
$u = w + q$,
\begin{align*}
  B_\ell^{(1)} &= 
  \partial\partial_{2}\Vb(w_{\ell}' + q'_{\ell - 1}, w_{\ell}' + q'_{\ell})
    \cdot
    \begin{pmatrix}
     u'_{\ell - 1} - w_{\ell}' - q_{\ell-1}' \\
     u'_{\ell} - w_{\ell}' - q'_{\ell}
    \end{pmatrix}
  \\ &\qquad +
  \partial\partial_{1}\Vb(w_{\ell}' + q'_{\ell}, w_{\ell}' + q'_{\ell + 1})
    \cdot
    \begin{pmatrix}
     u'_{\ell} - w_{\ell}' - q'_{\ell} \\
     u'_{\ell+1} - w_{\ell}' - q'_{\ell+1}
    \end{pmatrix}
  + \mathcal{O}(|w_{\ell - 1}''|^{2} + |w_{\ell}''|^{2})
  \\ &=
    \partial_{12}\Vb(w_{\ell}' + q_{\ell}', w_{\ell}' + q_{\ell + 1}')
       \cdot w_{\ell}''
       -
    \partial_{12}\Vb(w_{\ell}' + q_{\ell - 1}', w_{\ell}' + q_{\ell}') \cdot
w_{\ell - 1}''
  + \mathcal{O}(|w_{\ell - 1}''|^{2} + |w_{\ell}''|^{2}).
\end{align*}
Expanding the result again yields
\begin{displaymath}
  B_\ell^{(1)} =
 \bigl(\partial_{12}\Vb(w_{\ell}' + q_{\ell}', w_{\ell}' + q_{\ell}')
 \bigr) \cdot  (w_{\ell}'' - w_{\ell - 1}'')
 + \mathcal{O}\big(|q_{\ell}''| \cdot |w_{\ell - 1}''| +
               |q_{\ell}''| \cdot |w_{\ell}''| +
               |w_{\ell - 1}''|^{2} + |w_{\ell}''|^{2}\big).
\end{displaymath}
Thus, we obtain
\begin{equation}
  \label{eq:estiamte:B_ell^1}
  |B_\ell^{(1)}| \lesssim |w_{\ell}'''| + |w_{\ell - 1}''|^{2} + |w_{\ell}''|^{2}
   + |q_{\ell}''| \cdot |w_{\ell - 1}''| + |q_{\ell}''| \cdot |w_{\ell}''|.
\end{equation}

{\it Estimate for $B_\ell^{(2)} - B_\ell^{(3)}$: } The two terms
$B_\ell^{(2)}$ and $B_\ell^{(3)}$ are of similar structure and will need to be
treated together.  First, we rewrite $B_\ell^{(2)}$ and $B_\ell^{(3)}$ in the form
\begin{align*}
  B_\ell^{(2)}
  &=
    \int_{0}^{1}\partial \partial_{2}\Vb(w_{\ell}' + tq_{\ell - 1}',
                                          w_{\ell}' + tq_{\ell}')
   \cdot
   \begin{pmatrix}
     q_{\ell -1}' \\
     q_{\ell}'
   \end{pmatrix}dt
 \\ & \qquad +
 \int_{0}^{1}\partial \partial_{1}\Vb(w_{\ell}' + tq_{\ell}',
                                          w_{\ell}' + tq_{\ell + 1}')
   \cdot
   \begin{pmatrix}
     q_{\ell}' \\
     q_{\ell + 1}'
   \end{pmatrix}dt, \qquad \text{and} \\
  B_\ell^{(3)}
  &=
    \int_{0}^{1}\partial \partial_{2}\Vb(F_{0} + tq_{\ell - 1}',
                                          F_{0} + tq_{\ell}')
   \cdot
   \begin{pmatrix}
     q_{\ell -1}' \\
     q_{\ell}'
   \end{pmatrix}dt
 \\ &\qquad +
 \int_{0}^{1}\partial \partial_{1}\Vb(F_{0} + tq_{\ell}',
                                          F_{0} + tq_{\ell + 1}')
   \cdot
   \begin{pmatrix}
     q_{\ell}' \\
     q_{\ell + 1}'
   \end{pmatrix}dt.
\end{align*}
Subtracting and rearranging terms yields
\begin{align*}
  B_\ell^{(2)} - B_\ell^{(3)}  &=
    \int_{0}^{1}\biggl(\partial \partial_{2}\Vb(w_{\ell}' + tq_{\ell - 1}',
                                                 w_{\ell}' + tq_{\ell}')
                  - \partial \partial_{2}\Vb(F_{0} + tq_{\ell - 1}',
                                                 F_{0} + tq_{\ell}')\biggr)
   \cdot
   \begin{pmatrix}
     q_{\ell - 1}' \\
     q_{\ell}'
   \end{pmatrix}dt
 \\ &+
 \int_{0}^{1}\biggl(\partial \partial_{1} \Vb(w_{\ell}' + tq_{\ell - 1}',
                                                 w_{\ell}' + tq_{\ell}')
                  - \partial \partial_{1}\Vb(F_{0} + tq_{\ell}',
                                                 F_{0} + tq_{\ell + 1}')\biggr)
   \cdot
   \begin{pmatrix}
     q_{\ell}' \\
     q_{\ell + 1}'
   \end{pmatrix}dt.
\end{align*}
Making use of the Lipschitz continuity of $\partial^2 \Vb$ once again,
we deduce that
\begin{align*}
  \big| B_\ell^{(2)} - B_\ell^{(3)} \big|
  &\lesssim
 \int_{0}^{1}
   \begin{pmatrix}
    |w_{\ell}' - F_{0}| \\
    |w_{\ell}' - F_{0}|
   \end{pmatrix}
   \cdot
   \begin{pmatrix}
    q'_{\ell - 1} \\
    q'_{\ell}
   \end{pmatrix}
   dt
 +
 \int_{0}^{1}
    \begin{pmatrix}
    |w_{\ell}' - F_{0}| \\
    |w_{\ell}' - F_{0}|
   \end{pmatrix}
   \cdot
   \begin{pmatrix}
    q'_{\ell} \\
    q'_{\ell + 1}
   \end{pmatrix}
   dt
 \\ &\lesssim
      (|q_{\ell - 1}'| + |q_{\ell}'| + |q_{\ell + 1}'|)|w_{\ell}' - F_{0}|.
\end{align*}
Combining this estimate with \eqref{eq:estiamte:B_ell^1}, summing $\ell$, and
applying the Cauchy-Schwarz inequality yields the interior contribution of the
stated result.
\end{proof}

\begin{theorem}[Consistency]\label{Thm:PCConsistency}
  Let $\|f\|_{\U^*}$ be sufficiently small and $L$
  be sufficiently large so that the conditions of Theorems
  \ref{Thm:CBExistenceRegularity} and \ref{Thm:ExistenceInfiniteCorrector} with
  $F_0 = \nabla u^{\rm cb}(0)$ and of Proposition
  \ref{Prop:ExistenceFiniteCorrector} are satisfied. Let $u^{\rm cb}$ and $q_L$
  denote the corresponding Cauchy--Born and corrector solutions. Then,
  \begin{equation}
    \| \delta \E^{\rm a}(\Pi_{\rm a}u^{\rm cb} + q_{L}) - f\|_{\U^{*}}
    \lesssim
    \mu_{q}^{L}
    + |\nabla^2 u^{\rm cb}(0)|
    + \|\nabla^{2}u^{\rm cb}\|_{L^{4}}^{2} + \|\nabla^{3} u^{\rm cb}\|_{L^{2}}
    + \| \nabla f \|_{L^2}.
  \end{equation}
\end{theorem}
\begin{proof}
  As $u^{\rm cb}$ and $q_{\infty}$ are solutions to their respective problems,
  $\langle \delta \E^{\rm cb}(u^{\rm cb}), v \rangle
    = \langle \tilde{g}, \nabla v \rangle_{\R_{\geq 0}}$ for all
  $v \in \U^{\rm cb}$ with $\tilde{g}$ defined as in \eqref{eq:gtil=intf-cts}
  and $\langle \delta \E^{\Gamma}(q_{\infty}; F_0), v \rangle = 0$ for all
  $v \in \Q_{\infty}$.  We recall that $\U$ can be seen as a subspace of
  $\U^{\rm cb}$ and that $\Q_{\infty} = \U$.

  Let $\tilde{u} := \Pi_{\rm a}u^{\rm cb} + q_{L}$ and
  $F_0 = \nabla u^{\rm cb}(0)$, then using the fact that $u^{\rm cb}$ solves
  \eqref{SemiInfiniteCauchyBornProblem} we can split the consistency error into
  \begin{align*}
    \langle \delta \E^{\rm a}(\tilde{u}) - f, v \rangle
    &=
    \langle \delta_{\rm ext} + \delta_{\rm cb} + \delta_{\Gamma} + \delta_{\rm pc}, v \rangle, \qquad
    \text{where } \\
    \delta_{\rm ext}
    &:= \< f, \cdot \>_{\R_{\geq 0}} - \< f, \cdot \>_{\Z_{\geq 0}}, \\
    \delta_{\rm cb}
    &:=
      \delta \E^{\rm cb}(\Pi_{\rm a}u^{\rm cb}) - \delta \E^{\rm cb}(u^{\rm cb}), \\
    \delta_{\Gamma}
    &:=
      \delta \E^{\rm a}(\Pi_{\rm a}u^{\rm cb} + q_{L})
      - \delta \E^{\rm a}(\Pi_{\rm a}u^{\rm cb} + q_{\infty}), \qquad \text{and} \\
    \delta_{\rm pc}
    &:=
      \delta \E^{\rm a}(\Pi_{\rm a}u^{\rm cb} + q_{\infty})
      - \delta \E^{\rm cb}(\Pi_{\rm a}u^{\rm cb})
      - \delta \E^{\Gamma}(q_{\infty}; F_0).
  \end{align*}

  The term $\delta_{\rm ext}$ represents the error in the action of the
  external forces.  In the atomistic model, the external forces can be written
  as $\langle f, v \rangle_{\Z_\geq 0} = \langle g, v' \rangle_{\Z_{\geq 0}}$,
  where $g$ is defined as in \eqref{eq:g=sumf-discrete}.  Using
  \eqref{eq:estimate-g-gtil}, we get that
  \begin{equation}
    \< \delta_{\rm ext}, v \> =
    |\langle \tilde{g}, \nabla v \rangle_{\R_{\geq 0}}
     -  \langle g, v' \rangle_{\Z_{\geq 0}}|
   \lesssim \| \nabla f \|_{L^{2}} \| v' \|_{\ell^{2}}
   \qquad \text{ for all } v \in \U.
  \end{equation}

  The term $\delta_{\rm cb}$ represents the error associated with the
  discretization of the Cauchy--Born problem and was estimated in Theorem
  \ref{Proof:CBDiscretization} to be
  \begin{displaymath}
    \|\delta_{\rm cb}\|_{\U^*}
    \lesssim
     \|\nabla^{2}u^{\rm cb}\|_{L^{4}}^{2}.
  \end{displaymath}

  The term $\delta_{\Gamma}$ represents the error associated with the Galerkin
  projection of the corrector problem and was estimated in Theorem
  \ref{Prop:ExistenceFiniteCorrector}:
  \begin{displaymath}
    \| \delta_\Gamma \|_{\U^*} \lesssim \mu_q^L.
  \end{displaymath}

  Finally, $\delta_{\rm pc}$, the error due to the predictor-corrector method is
  estimated in Theorem \ref{Thm:SingleElementConsistencyError}:
  \begin{align*}
    \| \delta_{\rm pc} \|_{\U^*}
    &\lesssim
      \|(\Pi_{\rm a}u^{\rm
cb})''\|_{\ell^{4}}^{2}  +
        \|(\Pi_{\rm a}u^{\rm cb})'''\|_{\ell^{2}} +
        \|q'' \cdot (\Pi_{\rm a}u^{\rm cb})''\|_{\ell^{2}} \\
    & \qquad + \big|(\Pi_{\rm a}u^{\rm cb})'_0 - F_0 \big|
      + \Biggl(\sum_{\ell = 1}^{\infty}
          \bigl(|q_{\ell - 1}'| + |q_{\ell}'| + |q_{\ell + 1}'|\bigr)^{2}
          \cdot \big|(\Pi_{\rm a}u^{\rm cb})_{\ell}' -  F_0 \big|^{2}\Biggr)^{1/2} \\
    &=: P_1 + P_2 + P_3 + P_4 + P_5.
  \end{align*}
  where we have written $q \equiv q_\infty$ for the sake of simplicity of
  notation.

  To proceed, we need to relate discrete and continuous derivatives. For
  example,
  \begin{align*}
    \big| (\Pi_{\rm a} u^{\rm cb})_\ell'' \big|
    &= \big| (\Pi_{\rm a} u^{\rm cb})_{\ell+1}' - (\Pi_{\rm a} u^{\rm cb})_{\ell}' \big|
      = \bigg| \int_{\ell}^{\ell+1} \nabla u^{\rm cb}(x+1) - \nabla u^{\rm cb}(x)
      \, dx \bigg| \\
    & = \bigg| \int_{\ell}^{\ell+1} \int_0^1 \nabla^2 u^{\rm cb}(x+t) \,dt\,dx \bigg|
      \leq \int_{\ell}^{\ell+2} |\nabla^2 u^{\rm cb}(x)| \, dx.
  \end{align*}
  Analogously, we can prove
  \begin{align}
    \notag
    \big| (\Pi_{\rm a} u^{\rm cb})_\ell''' \big|
    &\leq \int_{\ell-1}^{\ell+2} |\nabla^3 u^{\rm cb}(x)| \, dx, \qquad \text{and} \\
    \label{eq:second-der-estimate-with-F0}
    \big|(\Pi_{\rm a}u^{\rm cb})'_\ell - F_0 \big|
    &\leq \int_0^{\ell+1} |\nabla^2 u^{\rm cb}(x)| \, dx.
  \end{align}
  Using these estimates, it is straightforward to establish that
  \begin{align*}
    P_1 = \| (\Pi_{\rm a} u^{\rm cb})'' \|_{\ell^4}^2
    \lesssim \| \nabla^2 u^{\rm cb} \|_{L^4}^2, \qquad \text{and}  \qquad
    P_2 = \| (\Pi_{\rm a} u^{\rm cb})''' \|_{\ell^2}
    \lesssim \| \nabla^3 u^{\rm cb} \|_{L^2}.
  \end{align*}

  The remaining three terms have similar structure. Estimating
  $|q''_\ell| \lesssim \mu_q^\ell = e^{- \alpha \ell}$ for some $\alpha > 0$, we
  deduce
  \begin{displaymath}
    P_3^2 = \|q'' \cdot (\Pi_{\rm a}u^{\rm cb})''\|_{\ell^{2}}^2 \lesssim
    \int_0^\infty e^{-2\alpha x} |\nabla^2 u^{\rm cb}(x)|^2 \, dx.
  \end{displaymath}
  Applying \eqref{eq:second-der-estimate-with-F0}, and using the fact that
  $e^{-\alpha x}$ is bounded below on $[0, 1]$, we have
  \begin{displaymath}
    P_4^2 = \big|(\Pi_{\rm a}u^{\rm cb})'_0 - F_0 \big|^2 \lesssim
    \int_0^1 e^{- 2\alpha x} |\nabla^2 u^{\rm cb}(x)|^2 \, dx.
  \end{displaymath}
  Finally, applying \eqref{eq:second-der-estimate-with-F0} again, together with
  $|q_\ell'| \lesssim \mu_q^\ell = e^{-\alpha \ell}$, we can estimate
  \begin{align*}
    P_5^2
    &= \sum_{\ell = 1}^{\infty}
      \bigl(|q_{\ell - 1}'| + |q_{\ell}'| + |q_{\ell + 1}'|\bigr)^{2}
      \cdot \big|(\Pi_{\rm a}u^{\rm cb})_{\ell}' -  F_0 \big|^{2} \\
    &\lesssim \sum_{\ell = 1}^\infty \mu_q^{2\ell} \, \bigg(\int_0^{\ell+1} |\nabla^2 u^{\rm cb}(x)| \, dx\bigg)^2
    \lesssim \int_0^\infty e^{-2\alpha t} \bigg( \int_0^t |\nabla^2 u^{\rm cb}(x)| \, dx \bigg)^2 dt \\
    &\lesssim \int_0^\infty e^{-2\alpha t} t \int_0^t |\nabla^2 u^{\rm cb}(x)|^2 \, dx \, dt
      = \int_0^\infty |\nabla^2 u^{\rm cb}(x)|^2 \int_{x}^\infty e^{-2\alpha t} t \, dt \, dx \\
    &\lesssim \int_0^\infty (1+x) e^{-2\alpha x} |\nabla^2 u^{\rm cb}(x)|^2 \, dx.
  \end{align*}
  Thus, we can combine
  \begin{displaymath}
    P_3^2 + P_4^2 + P_5^2 \lesssim \int_0^\infty (1+x) e^{-2\alpha x} |\nabla^2 u^{\rm cb}(x)|^2 \,dx.
  \end{displaymath}
  To finalize, we use the same argument as in the estimate for $P_5$ to deduce that
  \begin{align*}
    P_3^2 + P_4^2 + P_5^2
    &\lesssim |\nabla^2 u^{\rm cb}(0)|^2 +
      \int_0^\infty (1+x) e^{-2\alpha x} |\nabla^2 u^{\rm cb}(x) - \nabla^2 u^{\rm cb}(0)|^2 \, dx
      \\
    &\lesssim |\nabla^2 u^{\rm cb}(0)|^2 +
      \int_0^\infty (1+x)^2 e^{-2\alpha x} |\nabla^3 u^{\rm cb}(x)|^2 \, dx \\
    &\lesssim |\nabla^2 u^{\rm cb}(0)|^2 + \| \nabla^3 u^{\rm cb} \|_{L^2}^2.
  \end{align*}

  Combining the estimates for $\delta_{\rm ext}, \delta_{\rm cb}, \delta_\Gamma$
  and those for $P_j, j = 1, \dots, 5$, we arrive at the stated result.
\end{proof}

\begin{theorem}[Stability]\label{Thm:PCStability}
  There exists $\eps, L_0 > 0$ such that for all $f \in \U^{*}$ with
  $\|f\|_{\U^{*}} < \eps$ and all $L \geq L_{0}$, we have
  \begin{displaymath}
    \langle \delta^{2} \E^{\rm a}(\Pi_{\rm a}u^{\rm cb} + q_{L})v, v \rangle
    \geq \frac{c_{\rm a}}{2} \|v'\|_{\ell^{2}}^{2},
  \end{displaymath}
  where $c_{\rm a}$ denotes the stability constant for $u^{\rm a}_{\rm gr}$.
\end{theorem}
\begin{proof}
  First, let $f$ be small enough and $L$ be large enough so that Theorem
\ref{Thm:CBExistenceRegularity}, Theorem \ref{Thm:ExistenceInfiniteCorrector},
and Proposition \ref{Prop:ExistenceFiniteCorrector} holds.  The second
variation of the atomistic energy evaluated at the predictor-corrector solution
can be written as
 \begin{align*}
  \langle \delta^{2}& \E^{\rm a}(\Pi_{\rm a}u^{\rm cb} + q_{L})v, v \rangle
    =
    \langle \delta^{2} \E^{\rm a}(u^{\rm a}_{\rm gr})v, v \rangle
        - \bigl[\langle \delta^{2} \E^{\rm a}(u^{\rm a}_{\rm gr})v, v \rangle
         - \langle \delta^{2} \E^{\rm a}(q_{\infty})v, v \rangle\bigr]
    \\ &- \bigl[\langle \delta^{2} \E^{\rm a}(q_{\infty})v, v \rangle
         - \langle \delta^{2} \E^{\rm a}(q_{L})v, v \rangle\bigr]
        - \bigl[\langle \delta^{2} \E^{\rm a}(q_{L})v, v \rangle
         - \langle \delta^{2} \E^{\rm a}(\Pi_{\rm a}u^{\rm cb} + q_{L})v, v
\rangle\bigr].
 \end{align*}
 The highlighted difference terms can all be bound using the Lipschitz
continuity of the second variation of the atomistic energy:
 \begin{align}\label{Eq:StabilityTemp1}
   |\langle &\delta^{2} \E^{\rm a}(u^{\rm a}_{\rm gr})v, v \rangle
         - \langle \delta^{2} \E^{\rm a}(q_{\infty})v, v \rangle|
   \lesssim
   \|(u^{\rm a}_{\rm gr})' - (q_{\infty})'\|_{\ell^{2}}\|v'\|_{\ell^{2}}^{2}
   \\ \label{Eq:StabilityTemp2}
   |\langle &\delta^{2} \E^{\rm a}(q_{\infty})v, v \rangle
         - \langle \delta^{2} \E^{\rm a}(q_{L})v, v \rangle|
   \lesssim
   \|(q_{\infty})' - (q_{L})'\|_{\ell^{2}}\|v'\|_{\ell^{2}}^{2}
   \\ \label{Eq:StabilityTemp3}
   |\langle &\delta^{2} \E^{\rm a}(q_{L})v, v \rangle
         - \langle \delta^{2} \E^{\rm a}(\Pi_{\rm a}u^{\rm cb} + q_{L})v, v
\rangle|
   \lesssim
   \|\bigl(\Pi_{\rm a}u^{\rm cb}\bigr)'\|_{\ell^{2}}\|v'\|_{\ell^{2}}^{2}
 \end{align}
 Each of these bounds goes to zero either for large $L$ or for small
 $\|f\|_{\U^{*}}$.  The  bound in \eqref{Eq:StabilityTemp1} goes to zero as
 $\|f\|_{\U^{*}} \to 0$ by Theorem \ref{Thm:ExistenceInfiniteCorrector},
 the bound in \eqref{Eq:StabilityTemp2} goes to zero as $L \to \infty$ by
 Proposition \ref{Prop:ExistenceFiniteCorrector}, and the bound in
 \eqref{Eq:StabilityTemp3} goes to zero as $f \to 0$ by Theorem
 \ref{Thm:CBExistenceRegularity}.  Finally,
 $\langle \delta^{2} \E^{\rm a}(u^{\rm a}_{\rm gr})v, v \rangle \geq c_{\rm
   a}\|v'\|_{\ell^{2}}^{2}$
 by definition.  Thus, there exists an $L_{0} > 0$ and an $\eps > 0$ such that
 the theorem statement holds.
\end{proof}

\section{Appendix}

\begin{proposition}[First Variations.]\label{Prop:FirstVariations}  Let $u, v
\in \U$.  Then,
\begin{align*}
 \langle \delta \E^{\rm a}(u), v\rangle
   &=
 v_{0}'\biggl\{\partial_{1}V^{\surf}(u_{0}')
             + \partial_{1}\Vb(u_{0}',u_{1}')\biggr\}
 +
 \sum_{\ell = 1}^{\infty}v_{\ell}'
   \biggl\{\partial_{2}\Vb(u_{\ell - 1}',u_{\ell}')
         + \partial_{1}\Vb(u_{\ell}', u_{\ell + 1}')\biggr\},
 \\
 \langle \delta \E^{\Gamma}(q;u^{\rm cb}), v\rangle
   &=
  v_{0}'\biggl\{\partial_{1}V^{\surf}(F_{0} + q_{0}')
   + \partial_{1}\Vb(F_{0} + q_{0}', F_{0} + q_{1}')
   - \partial_{F}W(F_{0})\biggr\},
 \\ &+
 \sum_{\ell=1}^{\infty}v_{\ell}'
   \biggl\{\partial_{2}\Vb(F_{0} + q_{\ell - 1}', F_{0} + q_{\ell}')
         + \partial_{1}\Vb(F_{0} + q_{\ell}', F_{0} + q_{\ell + 1}')
 - \partial_{F}W(F_{0})\biggr\}.
\end{align*}
Let $u, v \in \U^{\rm cb}$.  Then,
\begin{displaymath}
 \langle \delta \E^{\rm cb}(u), v\rangle
   =
  \int_{0}^{\infty}\partial_{F}W(\nabla u) \cdot \nabla vdx
\end{displaymath}
\end{proposition}

\begin{proposition}[Second Variations.]\label{Prop:SecondVariations}
Let $u, v, w \in \U$.  Then,
\begin{align*}
 \langle \delta^{2} \E^{\rm a}(u)w, v \rangle
   &=
   v_{0}'F_{0}\Bigl\{\partial^{2}_{F}V^{\surf}(u_{0}') +
   \partial^{2}_{11}\Vb(u_{0}', u_{1}')\Bigr\}
 + v_{0}'w_{1}'\Bigl\{\partial^{2}_{12}\Vb(u_{0}',u_{1}')\Bigr\}
  \\
  + \sum_{\ell = 1}^{\infty}v_{\ell}'\biggl[
    &w_{\ell-1}'
      \Bigl\{\partial^{2}_{12}\Vb(u_{\ell-1}',u_{\ell}')\Bigr\}
   + w_{\ell}'\Bigl\{\partial^{2}_{22}\Vb(u_{\ell-1}',u_{\ell}') +
                     \partial^{2}_{11}\Vb(u_{\ell}',u_{\ell+1}')\Bigr\}
\\ + &w_{\ell+1}'
      \Bigl\{\partial^{2}_{12}\Vb(u_{\ell}',u_{\ell+1}')\Bigr\}\biggr]
\\
  \langle \delta^{2} \E^{\Gamma}(q;u^{\rm cb})w, v \rangle
   &=
   v_{0}'F_{0}\Bigl\{\partial^{2}_{F}V^{\surf}(F_{0} + q_{0}') +
   \partial^{2}_{11}\Vb(F_{0} + q_{0}', F_{0} + q_{1}')\Bigr\}
 \\ &+ v_{0}'w_{1}'\Bigl\{
     \partial^{2}_{12}\Vb(F_{0} + q_{0}', F_{0} + q_{1}')\Bigr\}
  \\
  + \sum_{\ell = 1}^{\infty}v_{\ell}'\biggl[
    &w_{\ell-1}'
      \Bigl\{
        \partial^{2}_{12}\Vb(F_{0} + q_{\ell-1}',F_{0} + q_{\ell}')
      \Bigr\}
\\ + &w_{\ell}'
     \Bigl\{
       \partial^{2}_{22}\Vb(F_{0} + q_{\ell-1}',F_{0}  q_{\ell}') +
       \partial^{2}_{11}\Vb(F_{0} + q_{\ell}',F_{0} + q_{\ell+1}')
     \Bigr\}
\\ + &w_{\ell+1}'
      \Bigl\{
        \partial^{2}_{12}\Vb(F_{0} + q_{\ell}',F_{0} + q_{\ell+1}')
      \Bigr\}\biggr]
\end{align*}
Let $u, v, w \in \U^{\rm cb}$.  Then,
\begin{displaymath}
 \langle \delta^{2} \E^{\rm cb}(u)w, v\rangle
   =
  \int_{0}^{\infty}\partial_{F}^{2}W(\nabla u) \cdot \nabla w \cdot \nabla vdx
\end{displaymath}
\end{proposition}

\subsection{Inverse Function Theorem}

\begin{theorem}\label{Thm:InverseFunctionTheorem}
  Let $\mathcal{H}$ be a Hilbert space equipped with norm
  $\| \cdot \|_{\mathcal{H}}$, and let
  $\mathcal{G} \in C^{1}(\mathcal{H},\mathcal{H}^{*})$ with Lipschitz-continuous
  derivative $\delta \mathcal{G}:$
 \begin{displaymath}
  \|\delta \mathcal{G}(u) - \delta \mathcal{G}(v)\|_{\mathcal{L}}
    \leq
  M\|u - v\|_{\mathcal{H}}
    \quad \text{for all} \;\; u,v \in \mathcal{H},
 \end{displaymath}
 where $\|\cdot\|_{\mathcal{L}}$ denotes the
 $\mathcal{L}(\mathcal{H},\mathcal{H}^{*})$-operator norm.

 Let $\bar{u} \in \mathcal{H}$ satisfy
 \begin{align*}
  \|\mathcal{G}(\bar{u})\|_{\mathcal{H}^{*}} &\leq \eta
    \\
  \langle \delta \mathcal{G}(\bar{u})v,v\rangle &\geq
    \gamma \| v\|^{2}_{\mathcal{H}}
    \quad \text{for all} \;\; v \in \mathcal{H},
 \end{align*}
 such that $M, \eta, \gamma$ satisfy the relation
 \begin{displaymath}
  \frac{2M\eta}{\gamma^{2}} < 1.
 \end{displaymath}
 Then, there exists a (locally unique) $u \in \mathcal{H}$
 such that $\mathcal{G}(u) = 0$,
 \begin{align*}
  \| u -  \bar{u}\|_{\mathcal{H}} &\leq 2\frac{\eta}{\gamma}, \quad \text{and}
    \\
  \langle \delta \mathcal{G}(u)v,v\rangle &\geq
    \left(1 - \frac{2M\eta}{\gamma^{2}}\right)\gamma\| v\|^{2}_{\mathcal{H}}
    \quad \text{for all} \; v \in \mathcal{H}.
 \end{align*}
\end{theorem}
\bibliography{PC_Method_References}{}

\begin{thebibliography}{10}

\bibitem{s00205-02-0218-5}
X.~Blanc, C.~Le~Bris, and P.-L. Lions.
\newblock From molecular models to continuum mechanics.
\newblock {\em Archive for Rational Mechanics and Analysis}, 164(4):341--381,
  2002.

\bibitem{PhysRevB.29.6443}
M.~S. Daw and M.~I. Baskes.
\newblock Embedded-atom method: Derivation and application to impurities,
  surfaces, and other defects in metals.
\newblock {\em Phys. Rev. B}, 29:6443--6453, Jun 1984.

\bibitem{NMAT:NMAT977}
J.~Diao, K.~Gall, and M.~L. Dunn.
\newblock Surface-stress-induced phase transformation in metal nanowires.
\newblock {\em Nature Materials}, 2(10):656--660, 2003.

\bibitem{m2an:2009007}
M.~Dobson and M.~Luskin.
\newblock An analysis of the effect of ghost force oscillation on
  quasicontinuum error.
\newblock {\em ESAIM: M2AN}, 43(3):591--604, 2009.

\bibitem{2010-JMPS-qce.stab}
M.~Dobson, M.~Luskin, and C.~Ortner.
\newblock Accuracy of quasicontinuum approximations near instabilities.
\newblock {\em J. Mech. Phys. Solids}, 58(10):1741--1757, 2010.

\bibitem{2010-MMS-SharpStabQCF}
M.~Dobson, M.~Luskin, and C.~Ortner.
\newblock Sharp stability estimates for the force-based quasicontinuum
  approximation of homogeneous tensile deformation.
\newblock {\em Multiscale Model. Simul.}, 8(3):782--802, 2010.

\bibitem{2010-ARMA-QCF}
M.~Dobson, M.~Luskin, and C.~Ortner.
\newblock Stability, instability, and error of the force-based quasicontinuum
  approximation.
\newblock {\em Arch. Ration. Mech. Anal.}, 197(1):179--202, 2010.

\bibitem{doi:10.1021/nl048456s}
K.~Gall, J.~Diao, and M.~L. Dunn.
\newblock The strength of gold nanowires.
\newblock {\em Nano Letters}, 4(12):2431--2436, 2004.

\bibitem{BF00261375}
M.~E. Gurtin and A.~I. Murdoch.
\newblock A continuum theory of elastic material surfaces.
\newblock {\em Archive for Rational Mechanics and Analysis}, 57(4):291--323,
  1975.

\bibitem{2013-ac.scb.1d}
K.~{Jayawardana}, C.~{Mordacq}, C.~{Ortner}, and H.~S. {Park}.
\newblock An analysis of the boundary layer in the {1D} surface {Cauchy-Born
  Model}.
\newblock {\em ESAIM: Math. Model. Numer. Anal.}, 47, 2013.

\bibitem{PhysRevB.37.3924}
R.~A. Johnson.
\newblock Analytic nearest-neighbor model for fcc metals.
\newblock {\em Phys. Rev. B}, 37:3924--3931, Mar 1988.

\bibitem{acta.atc}
M.~Luskin and C.~Ortner.
\newblock Atomistic-to-continuum coupling.
\newblock {\em Acta Numerica}, 22:397--508, 2013.

\bibitem{2011-MCOM-qnl1d}
C.~Ortner.
\newblock A priori and a posteriori analysis of the quasinonlocal
  quasicontinuum method in 1{D}.
\newblock {\em Math. Comp.}, 80(275):1265--1285, 2011.

\bibitem{2012-ARMA-cb}
C.~Ortner and F.~Theil.
\newblock Justification of the {Cauchy-Born} approximation of elastodynamics.
\newblock {\em Arch. Ration. Mech. Anal.}, 207, 2013.

\bibitem{PhysRevLett.95.255504}
H.~S. Park, K.~Gall, and J.~A. Zimmerman.
\newblock Shape memory and pseudoelasticity in metal nanowires.
\newblock {\em Phys. Rev. Lett.}, 95:255504, Dec 2005.

\bibitem{Park20083249}
H.~S. Park and P.~A. Klein.
\newblock A surface {Cauchy-Born} model for silicon nanostructures.
\newblock {\em Computer Methods in Applied Mechanics and Engineering},
  197(41-42):3249 -- 3260, 2008.
\newblock Recent Advances in Computational Study of Nanostructures.

\bibitem{NME:NME1754}
H.~S. Park, P.~A. Klein, and G.~J. Wagner.
\newblock A surface {Cauchy-Born} model for nanoscale materials.
\newblock {\em International Journal for Numerical Methods in Engineering},
  68(10):1072--1095, 2006.

\bibitem{TadmorMultiscaleBook}
E.~B. Tadmor and R.~E. Miller.
\newblock {\em Modeling Materials: Continuum, Atomistic, and Multiscale
  Techniques}.
\newblock Cambridge University Press, 2011.

\bibitem{EBTadmor:1996}
E.~B. Tadmor, M.~Ortiz, and R.~Phillips.
\newblock Quasicontinuum analysis of defects in solids.
\newblock {\em Philos. Mag. A}, 73:1529--1563, 1996.

\end{thebibliography}
\bibliographystyle{abbrv}

\end{document}